\documentclass[a4paper]{article}
\usepackage{amsfonts}
\usepackage{amsmath}
\usepackage{amssymb}
\begin{document}
\newtheorem{theorem}{Theorem}[section]
\newtheorem{lemma}{Lemma}[section]
\newtheorem{corollary}{Corollary}[section]
\newtheorem{conjecture}{Conjecture}[section]
\newtheorem{remark}{Remark}[section]
\newtheorem{definition}{Definition}[section]
\newtheorem{problem}{Problem}[section]
\newtheorem{example}{Example}[section]
\newtheorem{proposition}{Proposition}[section]
\title{{\bf PLURICANONICAL SYSTEMS OF PROJECTIVE VARIETIES OF GENERAL TYPE}}
\date{September, 2004}
\author{Hajime TSUJI}
\maketitle
\tableofcontents
\section{Introduction}

Let $X$ be a smooth projective variety and let $K_{X}$ be the canonical 
bundle of $X$. 
$X$ is said to be a general type., if there exists a positive 
integer $m$ such that the pluricanonical system 
$\mid mK_{X}\mid$ gives a birational (rational) embedding of $X$. 
The following problem is fundamental to study projective 
vareity of general type. \vspace{10mm}\\
{\bf Problem}
Let $X$ be a smooth projective variety of general type.
Find a positive integer $m_{0}$ such that for every $m\geqq m_{0}$,
$\mid mK_{X}\mid$ gives a birational rational map from $X$
into a projective space. \vspace{10mm} \\
If $\dim X = 1$, it is well known that $\mid 3K_{X}\mid$ gives a 
projective embedding.
In the case of smooth projective surfaces of general type, 
E. Bombieri showed that $\mid 5K_{X}\mid$ gives a birational
rational map from $X$ into a projective space (\cite{b3}).
But for the case of $\dim X\geqq 3$, very little is known about
the above problem.

The main purpose of this article is to prove the following theorems.

\begin{theorem} 
There exists a positive integer $\nu_{n}$ which depends
only on $n$ such that for every smooth projective $n$-fold $X$
of general type defined over complex numbers, $\mid mK_{X}\mid$ gives a birational rational map
from $X$ into a projective space for every $m\geqq \nu_{n}$.
\end{theorem}

Theorem 1.1 is very much related to the theory of minimal models.
It has been conjectured that for every nonuniruled smooth projective
variety $X$, there exists a projective variety $X_{min}$ such that 
\begin{enumerate}
\item $X_{min}$ is birationally equivalent to $X$,
\item $X_{min}$ has only {\bf Q}-factorial terminal singularities,
\item $K_{X_{min}}$ is a nef {\bf Q}-Cartier divisor.
\end{enumerate}
$X_{min}$ is called a minimal model of $X$. 
To construct a minimal model, the minimal model program 
(MMP) has been proposed (cf \cite[p.96]{k-m}).
The minimal model program was completed in the case of 
3-folds by S. Mori (\cite{mo}).

The proof of Theorem 1.1 can be very much simplified,
if we assume the existence of minimal models 
for projective varieties of general type.
The proof for the general case is modeled after the proof 
under the existence of minimal models  by using the theory of AZD. 

We should also note that even if we assume the existence
of minimal models for projective varieties of general type, 
Theorem 1.1 is quite nontrivial because the index of
a minimal model of $X$ 
(\cite[p.159, Definition 5.19]{k-m}) can be arbitrarily large. 
Conversely if we assume the existence of a minimal 
model of $X$ and  bound the index of the minimal model of 
$X$, then the proof of Theorem 1.1 is almost trivial. 
Since the index of minimal 3-folds of general type is 
unbounded, Theorem 1.1 is quite nontrivial even in the 
case of $\dim X = 3$. 
Hence in this sense the major difficulty of the proof of 
Theorem 1.1 is to find {\bf ``a (universal) lower bound'' of the positivity 
of} $K_{X}$.  
In fact Theorem 1.1 is equivalent to the following theorem. 
\begin{theorem}
There exists a positive number $C_{n}$ which depends
only on $n$ such that for every smooth projective $n$-fold $X$
of general type defined over complex numbers, 
\[
\mu (X,K_{X}) := n!\cdot\overline{\lim}_{m\rightarrow\infty}m^{-\dim X}\dim H^{0}(X,{\cal O}_{X}(mK_{X})) \geqq C_{n}
\]
holds.
\end{theorem} 
We note that $\mu (X,K_{X})$ is equal to the intersection 
number $K_{X}^{n}$ for a minimal projective $n$-fold
$X$ of general type (cf. Proposition 5.1 in Appendix). 
In Theorem 1.1 and 1.2, the numbers $\nu_{n}$ and 
$C_{n}$ have note yet been computed effectively 
except for the case of $n=3$ (\cite{3}). 

The relation of Theorem 1.1 and 1.2 is as follows.
Theorem 1.2 means that there exists a universal lower bound of 
the positivity of canonical bundle of smooth projective variety of 
general type with fixed dimension.
On the other hand, for a smooth projective variety of general type $X$, 
 the lower bound of $m$ such that $\mid mK_{X}\mid$ gives a birational embedding depends on the positivity of $K_{X}$ on certain families of subvarieties 
 which dominates $X$ as in Section 3@below. 
These families appears as the strata of the stratification as in 
\cite{t,a-s} which are the center of the log canonical 
singularities. 

The positivity of $K_{X}$ on the subvarieites  can be related 
to the positivity of the canonical bundles of the smooth model 
of the subvarities via the subadjunction theorem due to Kawamata (\cite{ka}). 
We note that since the family dominates $X$, for a general point of $X$, 
all the members of the family passing through the point should be of general type. 

The organization of the paper is as follows. 
In Section 2, we review the relation between multiplier ideal sheaves 
and singularity of divisors. 
In Section 3, we prove Theorem 1.1 and 1.2 assuming the existence 
of minimal models for projective varieties of general type. 
For the proof we use the induction on the dimension. 
Section 3.1 and 3.2 are similar to the argument as in \cite{t,a-s}. 
The essential part of Section 3 consists of Section 3.3 and 3.4.
In Section 3.3, we use the subadjunction theorem of Kawamata to relate 
the canonical divisor of centers of log canonical singularities and 
the canonical divisor of ambient space.
In Section 3.4, we prove that the minimal projective $n$-fold 
$X$ of general type with $K_{X}^{n}\leqq 1$ can be embedded
birationally into a projective space as a variety with degree $\leqq C^{n}$,
where $C$ is a positive constant depending only on $n$ (for the definition of 
$C$, see  Lemma 3.10).  Using this fact we finish the proof of 
Theorem 1.1 and 1.2 (assuming the existence of minimal models).

In Section 4, we prove Theorem 1.1 and 1.2 without assuming the 
existence of minimal models for projective varieties of general type.
Here we use the AZD (cf. Section 4.1) of $K_{X}$ instead of minimal models. 
The only essential difference of Section 3 and 4 is 
the use of the another subadjunction theorem (cf. Section 4.6,4.7,4.8). 
The rest is nothing but the transcription of Section 3 using 
the intersection theory of singular hermitian line bundles.  

In Section 5, we discuss the application of Theorem 1.1 and 1.2 
to Severi-Iitaka's conjecture. 

In this paper all the varieties are defined over {\bf C}.

\section{Multiplier ideal sheaves and singularities of divisors}

Before starting the proofs of Theorem 1.1 and 1.2, we shall review the 
relation between multiplier ideal sheaves and singularities of divisors.
In this subsection $L$ will denote a holomorphic line bundle on a complex manifold $M$. 
\begin{definition}
A  singular hermitian metric $h$ on $L$ is given by
\[
h = e^{-\varphi}\cdot h_{0},
\]
where $h_{0}$ is a $C^{\infty}$-hermitian metric on $L$ and 
$\varphi\in L^{1}_{loc}(M)$ is an arbitrary function on $M$.
We call $\varphi$ a  weight function of $h$.
\end{definition}
The curvature current $\Theta_{h}$ of the singular hermitian line
bundle $(L,h)$ is defined by
\[
\Theta_{h} := \Theta_{h_{0}} + \sqrt{-1}\partial\bar{\partial}\varphi ,
\]
where $\partial\bar{\partial}$ is taken in the sense of a current.
The $L^{2}$-sheaf ${\cal L}^{2}(L,h)$ of the singular hermitian
line bundle $(L,h)$ is defined by
\[
{\cal L}^{2}(L,h) := \{ \sigma\in\Gamma (U,{\cal O}_{M}(L))\mid 
\, h(\sigma ,\sigma )\in L^{1}_{loc}(U)\} ,
\]
where $U$ runs over the  open subsets of $M$.
In this case there exists an ideal sheaf ${\cal I}(h)$ such that
\[
{\cal L}^{2}(L,h) = {\cal O}_{M}(L)\otimes {\cal I}(h)
\]
holds.  We call ${\cal I}(h)$ the {\bf multiplier ideal sheaf} of $(L,h)$.
If we write $h$ as 
\[
h = e^{-\varphi}\cdot h_{0},
\]
where $h_{0}$ is a $C^{\infty}$ hermitian metric on $L$ and 
$\varphi\in L^{1}_{loc}(M)$ is the weight function, we see that
\[
{\cal I}(h) = {\cal L}^{2}({\cal O}_{M},e^{-\varphi})
\]
holds.
For $\varphi\in L^{1}_{loc}(M)$ we define the multiplier ideal sheaf of $\varphi$ by 
\[
{\cal I}(\varphi ) := {\cal L}^{2}({\cal O}_{M},e^{-\varphi}).
\]
\begin{example}
Let $\sigma\in \Gamma (X,{\cal O}_{X}(L))$ be the global section. 
Then 
\[
h := \frac{1}{\mid\sigma\mid^{2}} = \frac{h_{0}}{h_{0}(\sigma ,\sigma)}
\]
is a singular hemitian metric on $L$, 
where $h_{0}$ is an arbitrary $C^{\infty}$-hermitian metric on $L$
(the righthandside is ovbiously independent of $h_{0}$).
The curvature $\Theta_{h}$ is given by
\[
\Theta_{h} = 2\pi\sqrt{-1}(\sigma )
\]
where $(\sigma )$ denotes the current of integration over the 
divisor of $\sigma$. 
\end{example}
\begin{definition}
$L$ is said to be pseudoeffective, if there exists 
a singular hermitian metric $h$ on $L$ such that 
the curvature current 
$\Theta_{h}$ is a closed positive current.

Also a singular hermitian line bundle $(L,h)$ is said to be pseudoeffective, 
if the curvature current $\Theta_{h}$ is a closed positive current.
\end{definition}

If $\{\sigma_{i}\}$ are a finite number of global holomorphic sections of $L$,
for every positive rational number $\alpha$ and a $C^{\infty}$-function 
$\phi$,
\[
h := e^{-\phi}\cdot\frac{1}{\sum_{i}\mid\sigma_{i}\mid^{2\alpha}}
\]
defines a singular hermitian metric  on 
$\alpha L$.
We call such a metric $h$ a singular hermitian metric 
on $\alpha L$ with  {\bf algebraic singularities}.
Singular hermitian metrics with algebraic singularities 
are particulary easy to handle, because its multiplier 
ideal sheaf of the metric can 
be controlled by taking  suitable successive blowing ups 
such that the total transform of the divisor
$\sum_{i}(\sigma_{i})$ is a divisor with normal crossings. 

Let $D= \sum a_{i}D_{i}$ be an effective {\bf R}-divisor on $X$. 
Let $\sigma_{i}$ be a section of ${\cal O}_{X}(D_{i})$ with divisor $D_{i}$
respectively. 
Then we define 
\[
{\cal I}(D) : = {\cal I}(\prod_{i}\frac{1}{\mid\sigma_{i}\mid^{2a_{i}}})
\]
and call it the multiplier ideal of the divisor $D$.

Let us consider the relation between ${\cal I}(D)$ and 
singularities of $D$. 
\begin{definition}
Let $X$ be a normal variety and $D= \sum_{i}d_{i}D_{i}$ an effective {\bf Q}-divisor such that $K_{X}+D$ is {\bf Q}-Cartier. 
If $\mu : Y \longrightarrow X$ is an embedded resolution of the pair 
$(X,D)$, then we can write
\[
K_{Y} + \mu_{*}^{-1}D = \mu^{*}(K_{X}+D) + F
\]
with $F = \sum_{j}e_{j}E_{j}$ for the exceptional divisors $\{ E_{j}\}$.
We call $F$ the discrepancy and $e_{j}\in  \mbox{\bf Q}$ the discrepancy
coefficient for $E_{j}$. 
We regard $-d_{i}$ as the discrepancy coefficient of $D_{i}$. 

The pair $(X,D)$ is said to have only {\bf Kawamata log terminal singularities}
(KLT)(resp. {\bf log canonical singularities}(LC)), if 
$d_{i}< 1$(resp. $\leqq 1$) for all $i$ and $e_{j} > -1$ (resp. $\geqq -1$)
for all $j$ for an embedded resolution $\mu : Y \longrightarrow X$.
One can also say that $(X,D)$ is KLT (resp. LC), or $K_{X}+D$ is KLT
(resp. LC), when $(X,D)$ has only KLT (resp. LC).
The pair $(X,D)$ is said to be KLT (resp. LC) at a point $x_{0}\in X$,
if $(U,D\mid_{U})$ is KLT (resp. LC) for some neighbourhood $U$ 
of $x_{0}$. 
\end{definition}

The following proposition is a dictionary between algebraic geometry and 
the $L^{2}$-method.

\begin{proposition}
Let $D$ be a divisor on a smooth projective variety $X$. 
Then $(X,D)$ is KLT, if and only if 
${\cal I}(D)$ is trivial ($= {\cal O}_{X}$).
\end{proposition}
The proof is trivial and left to the readers. 
To locate the co-support of the multiplier ideal the following notion
is useful.  
\begin{definition}
A subvariety $W$ of $X$ is said to be a {\bf center of log canonical singularities} for the pair $(X,D)$, it there is a birational morphism from 
a normal variety $\mu : Y \longrightarrow X$ and a prime divisor $E$ on $Y$
with the discrepancy coefficient $e \leqq -1$
such that $\mu (E) = W$. 
\end{definition}
The set of all the centers of log canonical singularities is denoted 
by $CLC(X,D)$.
For a point $x_{0}\in X$, we define
$CLC(X,x_{0},D) := \{ W\in CLC(X,D)\mid x_{0}\in W\}$.
We quote the following proposition to introduce the notion of 
the minimal center of logcanoical singularities. 

\begin{proposition}(\cite[p.494, Proposition 1.5]{ka2})
Let $X$ be a  normal variety and $D$ an effective {\bf Q}-Cartier divisor 
such that $K_{X}+D$ is {\bf Q}-Cartier. 
Assume that $X$ is KLT and $(X,D)$ is LC.
If $W_{1},W_{2}\in CLC(X,D)$ and $W$ an irreducible component of $W_{1}\cap W_{2}$, then $W\in CLC(X,D)$. 
In particular if $(X,D)$ is not KLT at a point $x_{0}\in X$,
then there exists the unique minimal element of $CLC(X,x_{0},D)$. 
\end{proposition}

We call the minimal element  the {\bf minimal center of LC singularities} 
of $(X,D)$ at $x_{0}$.

\section{Proof of Theorem 1.1 and 1.2 assuming MMP}

In this section we prove Theorem 1.1 and 1.2 assuming the minimal model 
program.
The reason is that we can avoid inessential technicalities under 
this assumption.
The full proof of Theorem 1.1 and 1.2 is essentially nothing
but the transcription of the proof in this section 
by using the theory of AZD except the use of another 
subadjunction theorem (cf. Section 3.6, 3.8). 
Since the minimal model program is established in the case of 
3-folds, the proof under this assumption provides 
the full proofs of Theorem 1.1 and 1.2 for the case 
of projective varieties of general type of $\dim X \leqq 3$. 

Let $X$ be a minimal projective $n$-fold of general type, i.e.,
$X$ has only {\bf Q}-factorial terminal singularities and the 
canonical divisor $K_{X}$ is nef. 
We set 
\[
X^{\circ} = \{ x\in X_{reg} \mid x\not{\in} \mbox{Bs}\mid mK_{X}\mid \mbox{and  
$\Phi_{\mid mK_{X}\mid}$ is a biholomorphism} 
\]
\[
\hspace{50mm} \mbox{on a neighbourhood of $x$ for some $m\geqq 1$}\} .
\]
Then $X^{\circ}$ is a nonempty Zariski open subset of $X$. 
 
\subsection{Construction of a stratification}

The construction of a stratification below is similar to that 
in \cite{t,a-s}.  The only difference is the fact that we deal with
the {\bf Q}-Cartier divisor $K_{X}$ which is not Cartier in general. 
Of course this difference is very minor.

We set 
\[
\mu_{0}: = K_{X}^{n}.
\]
\begin{lemma} Let $x,x^{\prime}$ be distinct points on $X^{\circ}$.  
We set 
\[
{\cal M}_{x,x^{\prime}} = {\cal M}_{x}\otimes{\cal M}_{x^{\prime}},
\]
where ${\cal M}_{x},{\cal M}_{x^{\prime}}$ denote the
maximal ideal sheaves of the points $x,x^{\prime}$ respectively.
Let $\varepsilon$ be a sufficiently small positive number.
Then 
\[
H^{0}(X,{\cal O}_{X}(mK_{X})\otimes{\cal M}_{x,x^{\prime}}^{\lceil\sqrt[n]{\mu_{0}}
(1-\varepsilon )\frac{m}{\sqrt[n]{2}}\rceil})\neq 0
\]
for every sufficiently large $m$.
\end{lemma}
{\bf Proof of Lemma 3.1}.   
Let us consider the exact sequence:
\[
0\rightarrow H^{0}(X,{\cal O}_{X}(mK_{X})\otimes
{\cal M}_{x,x^{\prime}}^{\lceil\sqrt[n]{\mu_{0}}(1-\varepsilon )\frac{m}{\sqrt[n]{2}}\rceil})
\rightarrow H^{0}(X,{\cal O}_{X}(mK_{X}))\rightarrow
\]
\[
  H^{0}(X,{\cal O}_{X}
(mK_{X})\otimes {\cal O}_{X}/{\cal M}_{x,x^{\prime}}^{\lceil\sqrt[n]{\mu_{0}}(1-\varepsilon )\frac{m}{\sqrt[n]{2}}\rceil}).
\]
We note that 
\[
n!\cdot\overline{\lim}_{m\rightarrow\infty}m^{-n}\dim H^{0}(X,{\cal O}_{X}(mK_{X})) = \mu_{0}
\]
holds, since $K_{X}$ is nef and big (cf. Proposition 5.1 in 
Appendix). 
Since
\[
n!\cdot\overline{\lim}_{m\rightarrow\infty}m^{-n}\dim H^{0}(X,{\cal O}_{X}(mK_{X})
\otimes {\cal O}_{X}/{\cal M}_{x,x^{\prime}}^{\lceil\sqrt[n]{\mu_{0}}(1-\varepsilon )\frac{m}{\sqrt[n]{2}}\rceil})
=
\mu_{0}(1-\varepsilon )^{n} < \mu_{0}
\]
hold, we see that Lemma 3.1 holds.  {\bf Q.E.D.}

\vspace{5mm}

Let us take a sufficiently large positive integer $m_{0}$ and let $\sigma$
be a general (nonzero) element $\sigma_{0}$ of  
$H^{0}(X,{\cal O}_{X}(m_{0}K_{X})\otimes
{\cal M}_{x,x^{\prime}}^{\lceil\sqrt[n]{\mu_{0}}(1-\varepsilon )\frac{m_{0}}{\sqrt[n]{2}}\rceil})$. 
We define an effective {\bf Q}-divisor $D_{0}$ by 
\[
D_{0} = \frac{1}{m_{0}}(\sigma_{0}).
\]
We define a positive number $\alpha_{0}$ by
\[
\alpha_{0} := \inf\{\alpha > 0\mid 
\mbox{$(X,\alpha D_{0})$ is not KLT at 
both $x$ and $x^{\prime}$}\},
\]
where KLT is short for of Kawamata log terminal (cf. Definition 2.3).  
Since $(\sum_{i=1}^{n}\mid z_{i}\mid^{2})^{-n}$ is not locally integrable 
around $O\in \mbox{{\bf C}}^{n}$, by the construction of $D_{0}$, we see
that 
\[
\alpha_{0}\leqq \frac{n\sqrt[n]{2}}{\sqrt[n]{\mu_{0}}(1-\varepsilon )}
\]
holds.  About the relation between KLT condition and 
the multiplier ideal sheaves, please see Section 2.

Let us fix a positive number $\delta << 1$.
If we take $\varepsilon > 0$ sufficiently small, we may 
and do assume that 
\[
\alpha_{0} \leqq \frac{n\sqrt[n]{2}}{\sqrt[n]{\mu_{0}}}+ \delta 
\] 
holds. 
Then one of the following two cases occurs. \vspace{5mm} \\
{\bf Case} 1.1:  For every small positive number $\lambda$, 
$(X,(\alpha_{0}-\lambda )D_{0})$ is KLT 
at both $x$ and $x^{\prime}$. \\
{\bf Case} 1.2: For every small positive number $\delta$, 
$(X,(\alpha_{0}-\lambda )D_{0})$ is KLT  
at one of $x$ or $x^{\prime}$ say $x$. \vspace{10mm} \\

We first consider Case 1.1.
Let $X_{1}$ be the minimal center of log canonical singularities 
at $x$ (cf. Section 2). 
We consider the following two cases. \vspace{5mm} \\
{\bf Case} 3.1: $X_{1}$ passes through both $x$ and $x^{\prime}$, \\
{\bf Case} 3.2: Otherwise \vspace{10mm} \\

First we consider Case 3.1.
In this case $X_{1}$ is not isolated at $x$.  Let $n_{1}$ denote
the dimension of $X_{1}$.  Let us define the volume $\mu_{1}$ of $X_{1}$
with respect to $K_{X}$ by
\[
\mu_{1} := K_{X}^{n_{1}}\cdot X_{1}.
\]
Since $x\in X^{\circ}$, we see that $\mu_{1} > 0$ holds.
The proof of the following lemma is identical to that of Lemma 3.1.
\begin{lemma} Let $\varepsilon$ be a sufficiently small positive number and let $x_{1},x_{2}$ be distinct regular points on $X_{1}$. 
Then for a sufficiently large $m >1$,
\[
H^{0}(X_{1},{\cal O}_{X_{1}}(mK_{X})\otimes
{\cal M}_{x_{1},x_{2}}^{\lceil\sqrt[n_{1}]{\mu_{1}}(1-\varepsilon )\frac{m}{\sqrt[n_{1}]{2}}\rceil})\neq 0
\]
holds.
\end{lemma}
Let $x_{1},x_{2}$ be two distinct regular points on $X_{1}\cap X^{\circ}$. 
Let $m_{1}$ be a sufficiently large positive integer and 
Let 
\[
\sigma_{1}^{\prime}
\in 
H^{0}(X_{1},{\cal O}_{X_{1}}(m_{1}K_{X})\otimes
{\cal M}_{x_{1},x_{2}}^{\lceil\sqrt[n_{1}]{\mu_{1}}(1-\varepsilon )\frac{m}{\sqrt[n_{1}]{2}}\rceil})
\]
be a nonzero element. 

By Kodaira's lemma there is an effective {\bf Q}-divisor $E$ such
that $K_{X}- E$ is ample.
By the definition of $X^{\circ}$, we may assume that 
the support of $E$ does not contain both $x$ and $x^{\prime}$. 
Let $\ell_{1}$ be a sufficiently large positive integer which will be specified later  such that
\[
L_{1} := \ell_{1}(K_{X}- E)
\]
is Cartier. 

\begin{lemma}
If we take $\ell_{1}$ sufficiently large, then 
\[
\phi_{m} : H^{0}(X,{\cal O}_{X}(mK_{X}+L_{1}))\rightarrow 
H^{0}(X_{1},{\cal O}_{X_{1}}(mK_{X}+L_{1} ))
\]
is surjective for every  $m\geqq 0$.
\end{lemma}
{\bf Proof}.
$K_{X}$ is nef {\bf Q}-Cartier divisor by the assumption.
Let $r$ be the index of $X$, i.e. $r$ is the minimal positive integer such that $rK_{X}$ is Cartier.  Then for every locally free sheaf 
${\cal E}$, by Lemma 5.1 in Appendix,
there exists a positive integer $k_{0}$ 
depending on ${\cal E}$
such that for every $\ell \geqq k_{0}$ 
\[
H^{q}(X,{\cal O}_{X}((1+mr)K_{X}+L_{1})\otimes {\cal E}) = 0
\]
holds for every $q\geqq 1$ and $m\geqq 0$.
Let us consider the exact sequences 
\[
0 \rightarrow {\cal K}_{j}\rightarrow {\cal E}_{j}
\rightarrow {\cal O}_{X}(jK_{X})\otimes{\cal I}_{X_{1}}\rightarrow 0 
\]
for some locally free sheaf ${\cal E}_{j}$ for every 
$0 \leqq j \leqq r-1$,
where ${\cal I}_{X_{1}}$ denotes the ideal sheaf associated with
$X_{1}$. 
Then noting the above fact, we can prove that 
if  we take $\ell_{1}$ sufficiently large,
\[
H^{q}(X,{\cal O}_{X}(mK_{X}+L_{1})\otimes {\cal I}_{X_{1}}) = 0
\]
holds for every $q\geqq 1$ and $m\geqq 0$
by exactly the same manner as the standard proof of 
Serre's vanishing theorem (cf. \cite[p.228, Theorem 5.2]{ha}).
This implies the desired surjection. 
\vspace{5mm}{\bf Q.E.D.} \\ 

Let $\tau$ be a general section in 
$H^{0}(X,{\cal O}_{X}(L_{1}))$.
Then by Lemma 3.3 we see that   
\[
\sigma_{1}^{\prime}\otimes\tau\in
H^{0}(X_{1},{\cal O}_{X_{1}}(m_{1}K_{X}+L_{1})
{\cal M}_{x_{1},x_{2}}^{\lceil\sqrt[n_{1}]{\mu_{1}}(1-\varepsilon )\frac{m_{1}}
{\sqrt[n_{1}]{2}}\rceil})
\]
extends to a section
\[
\sigma_{1}\in H^{0}(X,{\cal O}_{X}((m_{1}+\ell_{1} )K_{X})).
\]
We may assume that  there exists a neighbourhood $U_{x,x^{\prime}}$ of $\{ x,x^{\prime}\}$ such that the divisor $(\sigma _{1})$  is smooth
on  $U_{x,x^{\prime}} - X_{1}$ by Bertini's theorem, if we take $\ell_{1}$
sufficiently large, since as in the proof of Lemma 3.3 
\[
H^{0}(X,{\cal O}_{X}(mK_{X}+L_{1}))
\rightarrow
H^{0}(X,{\cal O}_{X}(mK_{X}+L_{1})\otimes{\cal O}_{X}/{\cal I}_{X_{1}}\cdot{\cal M}_{y})
\]
is surjective for every $y\in X$ and
 $m\geqq 0$.
We set 
\[
D_{1} = \frac{1}{m_{1}+\ell_{1}}(\sigma_{1}). 
\]
Suppose that $x,x^{\prime}$ are {\bf nonsingular points} on $X_{1}$.
Then we set $x_{1} = x, x_{2} = x^{\prime}$.
Let $\varepsilon_{0}$ be a sufficiently small 
positive rational number and define $\alpha_{1}$ by
\[
\alpha_{1} := \inf\{\alpha > 0 \mid 
\mbox{ $(\alpha_{0}-\varepsilon_{0})D_{0} + \alpha D_{1}$ 
is not KLT at both $x$ and $x^{\prime}$} \}.
\]
Then we may define the proper subvariety $X_{2}$ of $X_{1}$ 
as a minimal center of log canonical singularities as before. 

\begin{lemma}
Let $\delta$ be the fixed positive number as above, then we may assume that 
\[
\alpha_{1}\leqq \frac{n_{1}\sqrt[n_{1}]{2}}{\sqrt[n_{1}]{\mu_{1}}} 
+ \delta 
\]
holds, if we make $\varepsilon_{0}$ and $\ell_{1}/m_{1}$ sufficiently 
small. 
\end{lemma}
To prove Lemma 3.4, we need the following elementary lemma.
\begin{lemma}(\cite[p.12, Lemma 6]{t})
Let $a,b$ be  positive numbers. Then
\[
\int_{0}^{1}\frac{r_{2}^{2n_{1}-1}}{(r_{1}^{2}+r_{2}^{2a})^{b}}
dr_{2}
=
r_{1}^{\frac{2n_{1}}{a}-2b}\int_{0}^{r_{1}^{-{2}{a}}}
\frac{r_{3}^{2n_{1}-1}}{(1 + r_{3}^{2a})^{b}}dr_{3}
\]
holds, where 
\[
r_{3} = r_{2}/r_{1}^{1/a}.
\]
\end{lemma}
{\bf Proof of Lemma 3.4.}
Let $(z_{1},\ldots ,z_{n})$ be a local coordinate on a 
neighbourhood $U$ of $x$ in $X$ such that 
\[
U \cap X_{1} = 
\{ q\in U\mid z_{n_{1}+1}(q) =\cdots = z_{n}(q)=0\} .
\] 
We set $r_{1} = (\sum_{i=n_{1}+1}^{n}\mid z_{1}\mid^{2})^{1/2}$ and 
$r_{2} = (\sum_{i=1}^{n_{1}}\mid z_{i}\mid^{2})^{1/2}$.
Fix an arbitrary $C^{\infty}$-hermitian metric $h_{X}$ on $K_{X}$. 
Then there exists a positive constant $C$ such that 
\[
\parallel\sigma_{1}\parallel^{2}\leqq 
C(r_{1}^{2}+r_{2}^{2\lceil\sqrt[n_{1}]{\mu_{1}}(1-\varepsilon )\frac{m_{1}}
{\sqrt[n_{1}]{2}}\rceil})
\]
holds on a neighbourhood of $x$, 
where $\parallel\,\,\,\,\parallel$ denotes the norm with 
respect to $h_{X}^{m_{1}+\ell_{1}}$.
We note that there exists a positive integer $M$ such that 
\[
\parallel\sigma_{1}\parallel^{-2} = O(r_{1}^{-M})
\]
holds on a neighbourhood of the generic point of $U\cap X_{1}$,
where $\parallel\,\,\,\,\parallel$ denotes the norm with respect to 
$h_{X}^{m_{0}}$. 
Then by Lemma 3.5, we have the inequality 
\[
\alpha_{1} \leqq (\frac{m_{1}+\ell_{1}}{m_{1}})\frac{n_{1}\sqrt[n_{1}]{2}}{\sqrt[n_{1}]{\mu_{1}}} + m_{1}\varepsilon_{0} 
\] 
holds.  Taking $\varepsilon_{0}$ and $\ell_{1}/m_{1}$ sufficiently small, 
we obtain that 
\[
\alpha_{1}\leqq \frac{n_{1}\sqrt[n_{1}]{2}}{\sqrt[n_{1}]{\mu_{1}}} 
+ \delta 
\]
holds.
{\bf Q.E.D.} \vspace{5mm} \\
If $x$ or $x^{\prime}$ is a singular point on $X_{1}$, we need the following lemma.
\begin{lemma}
Let $\varphi$ be a plurisubharmonic function on $\Delta^{n}\times{\Delta}$.
Let $\varphi_{t}(t\in\Delta )$ be the restriction of $\varphi$ on
$\Delta^{n}\times\{ t\}$.
Assume that $e^{-\varphi_{t}}$ does not belong to $L^{1}_{loc}(\Delta^{n},O)$
for every $t\in \Delta^{*}$.

Then $e^{-\varphi_{0}}$ is not locally integrable at $O\in\Delta^{n}$.
\end{lemma}
Lemma 3.6 is an immediate consequence of
the $L^{2}$-extension theorem \cite[p.20, Theorem]{o-t}.
Using Lemma 3.6 and Lemma 3.5, we see that Lemma 3.4 holds
by letting $x_{1}\rightarrow x$ and $x_{2}\rightarrow x^{\prime}$.

\vspace{5mm}

Next we consider Case 1.2 and Case 3.2.  
In this case  for every  sufficiently small positive number $\delta$, 
$(X,(\alpha_{0}-\varepsilon_{0})D_{0}+(\alpha_{1}-\delta )D_{1})$
is KLT at $x$ and not KLT at $x^{\prime}$. 

In these cases, instead of Lemma 3.2, we use the following simpler lemma.

\begin{lemma} Let $\varepsilon$ be a sufficiently small positive number and let $x_{1}$ be a smooth point on $X_{1}$. 
Then for a sufficiently large $m >1$,
\[
H^{0}(X_{1},{\cal O}_{X_{1}}(mK_{X})\otimes
{\cal M}_{x_{1}}^{\lceil\sqrt[n_{1}]{\mu_{1}}(1-\varepsilon )m\rceil})\neq 0
\]
holds.
\end{lemma}

Then taking a general nonzero element $\sigma_{1}^{\prime}$ in
\[
H^{0}(X_{1},{\cal O}_{X_{1}}(m_{1}K_{X})\otimes{\cal I}(h^{m_{1}})\otimes
{\cal M}_{x_{1}}^{\lceil\sqrt[n_{1}]{\mu_{1}}(1-\varepsilon )m_{1}
\rceil}),
\]
for a sufficiently large $m_{1}$.
As in Case 1.1 and Case 3.1 we obtain the proper subvariety
$X_{2}$ in $X_{1}$ also in this case.

Inductively for distinct points $x,x^{\prime}\in X^{\circ}$, we construct a strictly decreasing
sequence of subvarieties
\[
X = X_{0}\supset X_{1}\supset \cdots \supset X_{r}\supset X_{r+1} =  x\,\,\mbox{or}\,\, x^{\prime}
\]
and invariants (depending on small positive  rational numbers $\varepsilon_{0},\ldots ,
\varepsilon_{r-1}$, large positive integers $m_{0},m_{1},\ldots ,m_{r}$, etc.) :
\[
\alpha_{0} ,\alpha_{1},\ldots ,\alpha_{r},
\]
\[
\mu_{0},\mu_{1},\ldots ,\mu_{r}
\]
and
\[
n >  n_{1}> \cdots > n_{r}.
\]
By Nadel's vanishing theorem (\cite[p.561]{n}) we have the following lemma.
\begin{lemma} 
Let $x,x^{\prime}$ be two distinct points on $X^{\circ}$. 
Then for every $m\geqq \lceil\sum_{i=0}^{r}\alpha_{i}\rceil +1$,
$\Phi_{\mid mK_{X}\mid}$ separates $x$ and $x^{\prime}$.
And we may assume that 
\[
\alpha_{i} \leqq \frac{n_{i}\sqrt[n_{i}]{2}}{\sqrt[n_{i}]{\mu_{i}}}
+ \delta
\]
holds for every $0\leqq i\leqq r$. 
\end{lemma}
{\bf Proof}. 
For $i= 0,1,\ldots, r$ let $h_{i}$ be the singular hermitian metric 
on $K_{X}$ defined by 
\[
h_{i}:= \frac{1}{\mid\sigma_{i}\mid^{\frac{2}{m_{i}+\ell_{i}}}},
\]
where we have set $\ell_{0}:= 0$. 
More precisely for any $C^{\infty}$-hermitian metric $h_{X}$ 
on $K_{X}$ we have defined $h_{i}$ as 
\[
h_{i}:= \frac{h_{X}}{h_{X}^{m_{i}+\ell_{i}}(\sigma_{i},\sigma_{i})^{\frac{1}{m_{i}+\ell_{i}}}}.
\]
Using Kodaira's lemma (\cite[Appendix]{k-o}), 
let $E$ be an effective {\bf Q}-divisor $E$ such that
$K_{X} - E$ is ample. 
Let $m$ be a positive integer such that $m\geqq \lceil\sum_{i=0}^{r}\alpha_{i}\rceil +1$ holds. 
Let $h_{L}$ is a $C^{\infty}$-hermitian metric on the ample {\bf Q}-line bundle 
\[
L := (m-1-(\sum_{i=0}^{r-1}(\alpha_{i}-\varepsilon_{i}))
-\alpha_{r}-\delta_{L})K_{X} - \delta_{L}E
\]
 with strictly positive curvature, where $\delta_{L}$ be a sufficiently small positive number 
and we have considered $h_{L}$ as a singular hermitian metric 
on  \\ $(m-1- (\sum_{i=0}^{r-1}(\alpha_{i}-\varepsilon_{i}))-\alpha_{r})K_{X}$.
Let us define the singular hermitian metric $h_{x,x^{\prime}}$ of $(m-1)K_{X}$ defined by  
\[
h_{x,x^{\prime}} = (\prod_{i=0}^{r-1}h_{i}^{\alpha_{i}-\varepsilon_{i}})\cdot 
 h_{r}^{\alpha_{r}}\cdot h_{L}.
\]
Then we see that  ${\cal I}(h_{x,x^{\prime}})$ defines a subscheme of 
$X$ with isolated support around $x$ or $x^{\prime}$ by the definition of 
the invariants $\{\alpha_{i}\}$'s. 
By the construction the curvature current $\Theta_{h_{x,x^{\prime}}}$ is strictly positive on $X$. 
Then by Nadel's vanishing theorem (\cite[p.561]{n}) we see that 
\[
H^{1}(X,{\cal O}_{X}(mK_{X})\otimes {\cal I}(h_{x,x^{\prime}})) = 0
\]
holds. 
Hence 
\[
H^{0}(X,{\cal O}_{X}(mK_{X}))
\rightarrow 
H^{0}(X,{\cal O}_{X}(mK_{X})\otimes {\cal O}_{X}/{\cal I}(h_{x,x^{\prime}}))
\]
is surjective. 
Since by the construction of $h_{x,x^{\prime}}$ (if we take
$\delta_{L}$ sufficiently small) 
$\mbox{Supp}({\cal O}_{X}/{\cal I}(h_{x,x^{\prime}}))$ 
contains both $x$ and $x^{\prime}$ and is 
isolated at least one of $x$ or $x^{\prime}$. 
Hence by the above surjection, there exists a section
$\sigma\in H^{0}(X,{\cal O}_{X}(mK_{X}))$ such that 
\[
\sigma (x) \neq 0,\sigma (x^{\prime}) = 0
\]
or 
\[
\sigma (x) = 0,\sigma (x^{\prime}) \neq 0
\]
holds. 
This implies that $\Phi_{\mid mK_{X}\mid}$ separates 
$x$ and $x^{\prime}$.  

The proof of the last statement is similar to  the proof of Lemma 3.4.  {\bf Q.E.D.} 
\subsection{Construction of the stratification as a family}

In this subsection we shall construct the above stratification as a family. 
But this is not absolutely necessary for our proof of Theorem 1.1 and 1.2.
Please see Section 4.9 below for an alternative proof which bypasses this construction.

We note that for a fixed pair $(x,x^{\prime}) \in X^{\circ}\times X^{\circ}-\Delta_{X}$, $\sum_{i=0}^{r}\alpha_{i}$ depends on the choice of $\{ X_{i}\}$'s, where 
$\Delta_{X}$ denotes the diagonal of $X\times X$. 
Moving $(x,x^{\prime})$ in $X^{\circ}\times X^{\circ} - \Delta_{X}$, 
we  shall consider the above operation simultaneously.
Let us explain the procedure. 
We set 
\[
B := X^{\circ}\times X^{\circ} - \Delta_{X}.
\] 
Let 
\[
p : X\times B\longrightarrow X
\]
be the first projection and let   
\[
q : X\times {B}
\longrightarrow B
\]
be the second projection. 
Let $Z$ be the subvariety of $X\times B$ defined by
\[
Z := \{ (x_{1},x_{2},x_{3}) : X\times B \mid 
x_{1} = x_{2} \,\,\mbox{or} \,\, x_{1} = x_{3} \} .
\]
In this case we consider 
\[
q_{*}{\cal O}_{X\times B}(m_{0}p^{*}K_{X})
\otimes {\cal I}_{Z}^{\lceil\sqrt[n]{\mu_{0}}(1-\varepsilon )\frac{m_{0}}{\sqrt[n]{2}}\rceil }
\]
instead of 
\[
H^{0}(X,{\cal O}_{X}(m_{0}K_{X})\otimes
{\cal M}_{x,x^{\prime}}^{\lceil\sqrt[n]{\mu_{0}}(1-\varepsilon )\frac{m_{0}}{\sqrt[n]{2}}\rceil}),
\]
where ${\cal I}_{Z}$ denotes the ideal sheaf of $Z$. 
Let $\tilde{\sigma}_{0}$ be a nonzero global meromorphic 
section  of 
\[
q_{*}{\cal O}_{X\times B}(m_{0}p^{*}K_{X})
\otimes {\cal I}_{Z}^{\lceil\sqrt[n]{\mu_{0}}(1-\varepsilon )\frac{m_{0}}{\sqrt[n]{2}}\rceil } 
\]
on $B$ for a sufficiently large positive integer $m_{0}$.
We shall identify $\tilde{\sigma}_{0}$ with the family of sections 
of $m_{0}p^{*}K_{X}$.
We set 
\[
\tilde{D}_{0} : = \frac{1}{m_{0}}(\tilde{\sigma}_{0}).
\]
We define the singular hermitian metric $\tilde{h}_{0}$ 
on $p^{*}K_{X}$ by 
\[
\tilde{h}_{0}:= \frac{1}{\mid \tilde{\sigma}_{0}\mid^{2/m_{0}}}.
\]
We shall replace $\alpha_{0}$  by 
\[
\tilde{\alpha}_{0} 
:= \inf \{\alpha > 0\mid 
\mbox{the generic point of}\,\, Z \subseteq \mbox{Spec}
({\cal O}_{X \times B}/{\cal I}(h_{0}^{\alpha}))\} .
\]
Then for every $0 < \delta << 1$, there exists a Zariski 
open subset $U$ of $B$ such that for every $b \in U$, 
$\tilde{h}_{0}\mid_{X\times\{ b\}}$ is well defined and  
\[
b \not{\subseteq}\mbox{Spec}({\cal O}_{X\times\{ b\}}/{\cal I}(\tilde{h}_{0}^{\alpha_{0}-\delta}\mid_{X\times\{ b\}})),
\]
where we have identified $b$ with distinct two points in $X$. 
And also by Lemma 3.6, we see that 
\[
b \subseteq \mbox{Spec}({\cal O}_{X\times\{ b\}}/{\cal I}(\tilde{h}_{0}^{\alpha_{0}}\mid_{X\times\{ b\}})),
\]
holds for every $b\in B$. 
Let $\tilde{X}_{1}$ be the minimal center of log canonical singularities 
of $(X\times B,\alpha_{0}\tilde{D}_{0})$ at the generic point of $Z$. 
(although $\tilde{D}_{0}$ may not be effective this is meaningful 
by the construction of $\tilde{\sigma}_{0}$).   
We note that $\tilde{X}_{1}\cap q^{-1}(b)$ may not be 
irreducible even for a general $b\in B$. 
But if we take a suitable finite cover
\[
\phi_{0} : B_{0} \longrightarrow B,
\]
on the base change $X\times_{B}B_{0}$, $\tilde{X}_{1}$ 
defines a family of irreducible subvarieties
\[
f_{1} : \hat{X}_{1} \longrightarrow U_{0}
\]
of $X$ parametrized by a nonempty Zariski open subset
 $U_{0}$ of $\phi_{0}^{-1}(U)$.
Let $n_{1}$ be the relative dimension of $f_{1}$.
We set 
\[
\tilde{\mu}_{1} := K_{X}^{n_{1}}\cdot f_{1}^{-1}(b_{0})
\]
where $b_{0}$ is a general point on $U_{0}$.
Continuing this process 
we may construct a finite morphism 
\[
\phi_{r} : B_{r} \longrightarrow B
\]
and a nonempty Zariski open subset $U_{r}$ of $B_{r}$ 
which parametrizes a family of stratification 
\[
X \supset  X_{1} \supset X_{2} \supset \cdots \supset X_{r} 
\supset X_{r+1} =  x \,\,\mbox{or}\,\, x^{\prime}
\] 
constructed as before.
And we also obtain invariants $\{\tilde{\alpha}_{0},
\ldots ,\tilde{\alpha}_{r}\}$, $\{\tilde{\mu}_{0},\ldots ,\tilde{\mu}_{r}\}$,
$\{ n = \tilde{n}_{0}\ldots ,\tilde{n}_{r}\}$.
Hereafter we denote these invariants without $\,\,\tilde{} \,\,$ for simplicity.  By the same proof as in  Lemma 3.4, we have the following lemma. 
\begin{lemma}
We may assume that 
\[
\alpha_{i}\leqq \frac{n_{i}\sqrt[n_{i}]{2}}{\sqrt[n_{i}]{\mu_{i}}} + \delta
\]
holds for every $0 \leqq i\leqq r$.
\end{lemma}
By Lemma 3.8, we obtain that 
For every 
\[
m > \lceil\sum_{i=0}^{r}\alpha_{i}\rceil + 1
\]
$\mid mK_{X}\mid$ separates points on the nonempty Zariski 
open subset $\phi_{r}(B_{r})$.

Let us consider the complement 
$X^{\circ} - \phi_{r}(B_{r})$.
Replacing $B$ by 
\[
(X^{\circ} - \phi_{r}(B_{r}))\times 
(X^{\circ} - \phi_{r}(B_{r})) - \Delta_{X},
\] 
we may continue the 
same process. 
Hence by Noetherian induction, we have 
the following proposition.
\begin{proposition}
There exists a finite 
stratification of $X \times X - \Delta_{X}$ such that 
the each stratum supports a (multivalued) family of stratification 
of $X$ : 
\[
X \supset  X_{1} \supset X_{2} \supset \cdots \supset X_{r} 
\supset X_{r+1} =  x \,\,\mbox{or}\,\, x^{\prime}
\] 
with the same invariants $\{ \alpha_{0},\cdots ,\alpha_{r}\}$,
$\{ \mu_{0},\cdots ,\mu_{r}\}$ etc. ($r$ may depend on the strata).
\end{proposition}
\subsection{Use of Kawamata's subadjunction theorem}
The following subadjunction theorem is crucial in our proof. 
\begin{theorem}(\cite{ka})
Let $X$ be a normal projective variety.
Let $D^{\circ}$ and $D$ be effective {\bf Q}-divisor on $X$ such that 
$D^{\circ} < D$, $(X,D^{\circ})$ is log terminal and 
$(X,D)$ is log canonical. 
Let $W$ be a minimal center of log canonical singularities for $(X,D)$. 
Let $H$ be an ample Cartier divisor on $X$ and $\epsilon$ a positive rational number.
Then there exists an effective {\bf Q}-divisor $D_{W}$ on $W$ such that 
\[
(K_{X}+D+\epsilon H)\mid_{W}\sim_{\mbox{\bf Q}}K_{W}+D_{W}
\]
and $(W,D_{W})$ is log terminal. 
In particular $W$ has only rational singularities.
\end{theorem}
\begin{remark}
As is stated in \cite[Remark 3.2]{ka3}, 
the assumption that $W$ is a minimal center can be replaced
that $W$ is a local minimal center, since the argument 
in \cite{ka} which uses the variation of Hodge structure 
does not change. 
But in this case we need to replace $K_{W}$ by the 
pushforward of the canonical divisor of the normalization 
of $W$, 
since  $W$ may be nonnormal (\cite[p.494, Theorem 1.6]{ka2} 
works only locally in this case). 
\end{remark}
Roughly speaking, Theorem 3.1 implies that $K_{X}+D\mid_{W}$ (almost) dominates $K_{W}$.  

Let us consider again the sequence of numbers $\alpha_{j}$, divisors $D_{j}$
and the stratification $X \supset X_{1} \supset \cdots \supset X_{j}$
 which were defined in Section 3.1.
Let $W_{j}$ be a nonsingular model of $X_{j}$. 
Applying Theorem 3.1 to $K_{X}+D$ where
\[
D = (\alpha_{0}-\varepsilon_{0})D_{0} +\cdots +
(\alpha_{j-2}-\varepsilon_{j-2})D_{j-2} + \alpha_{j-1}D_{j-1},
\]
we get 
\[
\mu (W_{j},K_{W_{j}}) \leqq \mu (W_{j},K_{W_{j}}+D_{W_{j}}) 
\leqq 
(1 +\sum_{i=0}^{j-1}\alpha_{i})^{n_{j}}\cdot \mu_{j}
\]
hold, where 
\[
\mu (W_{j},K_{W_{j}}) 
:= n_{j}!\cdot \overline{\lim}_{m\rightarrow\infty}
m^{-n_{j}}\dim H^{0}(W_{j},{\cal O}_{W_{j}}(mK_{W_{j}})).
\]
We note that if we take $x,x^{\prime}$ general, 
$W_{j}$ ought to be of general type. 
More precisely there exists no subfamily ${\cal K}$ of 
$\{ W_{j}\}$ (parametrized a quasiprojective variety) such that 
\begin{enumerate}
\item every member of ${\cal K}$ is of non-general type, 
\item the members of ${\cal K}$ dominates $X$ by the natural morphism. 
\end{enumerate}
Otherwise $X$ is dominated by a family of varieties 
of nongeneral type and 
this contradicts the assumption that $X$ is of general type. 
Hence there exists a nonempty Zariski open set $U_{0}$ of $X^{\circ}$ 
such that if $(x,x^{\prime}) \in U_{0}\times U_{0}$, then  $W_{j}$ is of general
type for every $j$.

We shall prove Theorem 1.2 by  induction on $n$. 
Suppose that Theorem 1.2 holds for projective varieties of 
general type of dimension less than or equal to $n-1$ (the case of $n=1$ is trivial), i.e.,
for every positive integer $k < n$
there exists a positive number $C(k)$ such that 
for every smooth projective variety $W$ of general type 
of dimension $k$, 
\[
\mu (W,K_{W})\geqq C(k)
\]
holds. 

Let us consider again the sequence of numbers $\alpha_{j}$, divisors $D_{j}$
and the stratification $X \supset X_{1} \supset \cdots \supset X_{j}$
 which were defined in Section 3.1.
Then by the above inequality  
\[
C(n_{j}) \leqq 
(1 +\sum_{i=0}^{j-1}\alpha_{i})^{n_{j}}\cdot \mu_{j}
\]
holds. 
Since 
\[
\alpha_{i} \leqq \frac{\sqrt[n_{i}]{2}\, n_{i}}{\sqrt[n_{i}]{\mu_{i}}}+ \delta
\]
holds by Lemma 3.9, 
we see that  
\[
\frac{1}{\sqrt[n_{j}]{\mu_{j}}}\leqq (1+\sum_{i=0}^{j-1}\frac{\sqrt[n_{i}]{2}\, n_{i}}{\sqrt[n_{i}]{\mu_{i}}})\cdot C(n_{j})^{-\frac{1}{n_{j}}}
\]
holds for every $j \geqq 1$.
We recall the finite stratification of $X^{\circ}\times X^{\circ}
-\Delta_{X}$ in Section 3.2. 
Using the above inequality inductively, we have
the following lemma.
\begin{lemma}
Suppose that $\mu_{0} \leqq 1$ holds.
Then there exists a positive constant $C$ depending only on $n$ 
such that for every $(x,x^{\prime})\in U_{0}\times U_{0} -\Delta_{X}$ 
the corresponding invariants $\{ \mu_{0},\cdots ,\mu_{r}\}$ 
and $\{ n_{1},\cdots ,n_{r}\}$  depending on $(x,x^{\prime})$
($r$ may also depend on $(x,x^{\prime})$) satisfies
the inequality :
\[
1+\lceil \sum_{i=0}^{r}\frac{\sqrt[n_{i}]{2}\, n_{i}}{\sqrt[n_{i}]{\mu_{i}}}\rceil \leqq \lfloor\frac{C}{\sqrt[n]{\mu_{0}}}\rfloor .
\]
\end{lemma} 
\subsection{Estimate of the degree}
To relate $\mu_{0}$ and the degree of the pluricanonical image of $X$, 
we need the following lemma. 
\begin{lemma}
If $\Phi_{\mid mK_{X}\mid}\mid$ is birational rational map
onto its image, then
\[
\deg \Phi_{\mid mK_{X}\mid}(X)\leqq \mu_{0}\cdot m^{n}
\]
holds.
\end{lemma}
{\bf Proof}.
Let $p : \tilde{X}\longrightarrow X$ be the resolution of 
the base locus of $\mid mK_{X}\mid$ and let 
\[
p^{*}\mid mK_{X}\mid = \mid P_{m}\mid + F_{m}
\]
be the decomposition into the free part $\mid P_{m}\mid$ 
and the fixed component $F_{m}$. 
We have
\[
\deg \Phi_{\mid mK_{X}\mid}(X) = P_{m,}^{n}
\]
holds.
Then by the ring structure of $R(X,K_{X})$, we have an injection 
\[
H^{0}(\tilde{X},{\cal O}_{\tilde{X}}(\nu P_{m}))\rightarrow 
H^{0}(X,{\cal O}_{X}(m\nu K_{X}))
\]
for every $\nu\geqq 1$.
We note that since ${\cal O}_{\tilde{X}}(\nu P_{m})$ is globally generated
on $\tilde{X}$, for every $\nu \geqq 1$ we have the injection 
\[
{\cal O}_{\tilde{X}}(\nu P_{m})\rightarrow p^{*}{\cal O}_{X}(m\nu K_{X}).
\]
Hence there exists a natural morphism 
\[
H^{0}(\tilde{X},{\cal O}_{\tilde{X}}(\nu P_{m}))
\rightarrow 
H^{0}(X,{\cal O}_{X}(m\nu K_{X}))
\]
for every $\nu\geqq 1$. 
This morphism is clearly injective. 
This implies that 
\[
\mu_{0} \geqq  m^{-n}\mu (\tilde{X}_{i},P_{m})
\]
holds. 
Since $P_{m}$ is nef and big on $X$, we see that 
\[
\mu (\tilde{X},P_{m}) = P_{m}^{n}
\]
holds.
Hence
\[
\mu_{0}\geqq m^{-n}P_{m}^{n}
\]
holds.  This implies that
\[
\deg \Phi_{\mid mK_{X}\mid}(X)\leqq \mu_{0}\cdot m^{n}
\]
holds.
{\bf Q.E.D.}
\vspace{5mm} \\
\subsection{Completion of the proof of Theorem 1.1 and 1.2 assuming MMP}
By Lemma 3.9.3.10 and 3.11 we see that 
if $\mu_{0} \leqq 1$ holds, 
for 
\[
m := \lfloor\frac{C}{\sqrt[n]{\mu_{0}}}\rfloor ,
\]
$\mid mK_{X}\mid$ gives a birational embedding of $X$ and 
\[
\deg \Phi_{\mid mK_{X}\mid}(X) 
\leqq C^{n}
\]
holds, where $C$ is the positive constant in Lemma 3.10.  
Also 
\[
\dim H^{0}(X,{\cal O}_{X}(mK_{X}))
\leqq n+1 + \deg \Phi_{\mid mK_{X}\mid}(X) 
\]
holds by the semipositivity of the $\Delta$-genus (\cite{fu}). 
Hence we have that if $\mu_{0} \leqq 1$,
\[
\dim H^{0}(X,{\cal O}_{X}(mK_{X})) 
\leqq n+1+ C^{n} 
\]
holds. 

Since $C$ is a positive constant depending only on $n$, 
combining  the above two inequalities, 
we have that there exists a positive constant $C(n)$ 
depending only on $n$ such that 
\[
\mu_{0} = K_{X}^{n} \geqq C(n)
\]
holds. 

More precisely we argue as follows. 
Let ${\cal H}$ be an irreducible component of the Hilbert scheme of a projective spaces of dimension $\leqq n+ C^{n}$.  
Let ${\cal H}_{0}$ be the Zariski open subset of ${\cal H}$ which parametrizes 
irreducible subvarieties. 
Then there exists a finite stratification of ${\cal H}_{0}$ by Zariski locally closed subsets such that on each stratum there exists a simultaneous resolution 
of the universal family on the strata. 
We note that the volume of the canonical bundle of the resolution is constant on each strata by \cite{tu6,nak}.
Hence there exists a positive constant $C(n)$ depending only on $n$  such that 
\[
\mu (X,K_{X})\geqq C(n)
\]
holds for every projective $n$-fold $X$ of general type 
by the degree bound as above in the case of $\mu_{0}\leqq 1$. 
This completes the proof of Theorem 1.2
assuming MMP.

Then by Lemma 3.9 and 3.10, we see that there exists 
a positive integer $\nu_{n}$ depending only on $n$ such that 
for every projective $n$-fold $X$ of general type, 
$\mid mK_{X}\mid$ gives a birational embedding into a 
projective space for every $m\geqq \nu_{n}$. 
This completes the proof of Theorem 1.1 assuming MMP.

\section{Proof of Theorem 1.1 and 1.2 without assuming 
MMP}

In this section we shall prove Theorem 1.1 and 1.2 in full
generality. The proof is almost parallel to the one  assuming 
MMP, if we replace the minimal model by an AZD (analytic
Zariski decomposition) of the canonical line bundle. 

\subsection{Analytic Zariski decomposition}

To study a pseudoeffective line bundle we introduce the notion of analytic Zariski
decompositions.
By using analytic Zariski decompositions, we can handle a pseudoeffective line bundle, as if it  were a nef line bundle.
\begin{definition}
Let $M$ be a compact complex manifold and let $L$ be a line bundle
on $M$.  A singular hermitian metric $h$ on $L$ is said to be 
an analytic Zariski decomposition, if the followings hold.
\begin{enumerate}
\item $\Theta_{h}$ is a closed positive current,
\item for every $m\geqq 0$, the natural inclusion
\[
H^{0}(M,{\cal O}_{M}(mL)\otimes{\cal I}(h^{m}))\rightarrow
H^{0}(M,{\cal O}_{M}(mL))
\]
is isomorphim.
\end{enumerate}
\end{definition}
\begin{remark} If an AZD exists on a line bundle $L$ on a smooth projective
variety $M$, $L$ is pseudoeffective by the condition 1 above.
\end{remark}

\begin{theorem}(\cite{tu,tu2})
 Let $L$ be a big line  bundle on a smooth projective variety
$M$.  Then $L$ has an AZD. 
\end{theorem}
As for the existence for general pseudoeffective line bundles, 
now we have the following theorem.
\begin{theorem}(\cite[Theorem 1.5]{d-p-s})
Let $X$ be a smooth projective variety and let $L$ be a pseudoeffective 
line bundle on $X$.  Then $L$ has an AZD.
\end{theorem}
Although the proof is in \cite{d-p-s}, 
we shall give a proof here, because we shall use it afterward. 

 Let  $h_{0}$ be a fixed $C^{\infty}$-hermitian metric on $L$.
Let $E$ be the set of singular hermitian metric on $L$ defined by
\[
E = \{ h ; h : \mbox{lowersemicontinuous singular hermitian metric on $L$}, 
\]
\[
\hspace{70mm}\Theta_{h}\,
\mbox{is positive}, \frac{h}{h_{0}}\geq 1 \}.
\]
Since $L$ is pseudoeffective, $E$ is nonempty.
We set 
\[
h_{L} = h_{0}\cdot\inf_{h\in E}\frac{h}{h_{0}},
\]
where the infimum is taken pointwise. 
The supremum of a family of plurisubharmonic functions 
uniformly bounded from above is known to be again plurisubharmonic, 
if we modify the supremum on a set of measure $0$(i.e., if we take the uppersemicontinuous envelope) by the following theorem of P. Lelong.

\begin{theorem}(\cite[p.26, Theorem 5]{l})
Let $\{\varphi_{t}\}_{t\in T}$ be a family of plurisubharmonic functions  
on a domain $\Omega$ 
which is uniformly bounded from above on every compact subset of $\Omega$.
Then $\psi = \sup_{t\in T}\varphi_{t}$ has a minimum 
uppersemicontinuous majorant $\psi^{*}$  which is plurisubharmonic.
We call $\psi^{*}$ the uppersemicontinuous envelope of $\psi$. 
\end{theorem}
\begin{remark} In the above theorem the equality 
$\psi = \psi^{*}$ holds outside of a set of measure $0$(cf.\cite[p.29]{l}). 
\end{remark}

By Theorem 4.3,we see that $h_{L}$ is also a 
singular hermitian metric on $L$ with $\Theta_{h}\geq 0$.
Suppose that there exists a nontrivial section 
$\sigma\in \Gamma (X,{\cal O}_{X}(mL))$ for some $m$ (otherwise the 
second condition in Definition 4.1 is empty).
We note that  
\[
\frac{1}{\mid\sigma\mid^{\frac{2}{m}}} 
\]
gives the weight of a singular hermitian metric on $L$ with curvature 
$2\pi m^{-1}(\sigma )$, where $(\sigma )$ is the current of integration
along the zero set of $\sigma$. 
By the construction we see that there exists a positive constant 
$c$ such that  
\[
\mbox{($\star$)}\hspace{10mm} \frac{h_{0}}{\mid\sigma\mid^{\frac{2}{m}}} \geq c\cdot h_{L}
\]
holds. 
Hence
\[
\sigma \in H^{0}(X,{\cal O}_{X}(mL)\otimes{\cal I}_{\infty}(h_{L}^{m}))
\]
holds.  
Hence in praticular
\[
\sigma \in H^{0}(X,{\cal O}_{X}(mL)\otimes{\cal I}(h_{L}^{m}))
\]
holds.  
 This means that $h_{L}$ is an AZD of $L$. 
\vspace{10mm} {\bf Q.E.D.} 

\begin{remark}
By the above proof (cf. ($\star$)) we have that for the AZD $h_{L}$ constructed 
as above
\[
H^{0}(X,{\cal O}_{X}(mL)\otimes{\cal I}_{\infty}(h_{L}^{m}))
\simeq 
H^{0}(X,{\cal O}_{X}(mL))
\]
holds for every $m$, where ${\cal I}_{\infty}(h_{L}^{m})$ 
denotes the $L^{\infty}$-multiplier ideal sheaf, i.e., 
for every open subset $U$ in $X$, 
\[
{\cal I}_{\infty}(h_{L}^{m})(U):= \{ f\in {\cal O}_{X}(U) \mid
\mid f\mid^{2}(h_{L}/h_{0})^{m}\in L^{\infty}_{loc}(U)\} .
\]
\end{remark}

\subsection{The $L^{2}$-extension theorem}

 Let $M$ be a complex manifold of dimension $n$ and let $S$ be a closed complex submanifold of $M$. 
Then we consider a class of continuous function $\Psi : M\longrightarrow [-\infty , 0)$  such that  
\begin{enumerate}
\item $\Psi^{-1}(-\infty ) \supset S$,
\item if $S$ is $k$-dimensional around a point $x$, there exists a local 
coordinate $(z_{1},\ldots ,z_{n})$ on a neighbourhood of $x$ such that 
$z_{k+1} = \cdots = z_{n} = 0$ on $S\cap U$ and 
\[
\sup_{U\backslash S}\mid \Psi (z)-(n-k)\log\sum_{j=k+1}^{n}\mid z_{j}\mid^{2}\mid < \infty .
\]
\end{enumerate} 
The set of such functions $\Psi$ will be denoted by $\sharp (S)$. 

For each $\Psi \in \sharp (S)$, one can associate a positive measure 
$dV_{M}[\Psi ]$ on $S$ as the minimum element of the 
partially ordered set of positive measures $d\mu$ 
satisfying 
\[
\int_{S_{k}}f\, d\mu \geqq 
\overline{\lim}_{t\rightarrow\infty}\frac{2(n-k)}{v_{2n-2k-1}}
\int_{M}f\cdot e^{-\Psi}\cdot \chi_{R(\Psi ,t)}dV_{M}
\]
for any nonnegative continuous function $f$ with 
$\mbox{supp}\, f\subset\subset M$.
Here $S_{k}$ denotes the $k$-dimensional component of $S$,
$v_{m}$ denotes the volume of the unit sphere 
in $\mbox{\bf R}^{m+1}$ and 
$\chi_{R(\Psi ,t)}$ denotes the characteristic funciton of the set 
\[
R(\Psi ,t) = \{ x\in M\mid -t-1 < \Psi (x) < -t\} .
\]

Let $M$ be a complex manifold and let $(E,h_{E})$ be a holomorphic hermitian vector 
bundle over $M$. 
Given a positive measure $d\mu_{M}$ on $M$,
we shall denote $A^{2}(M,E,h_{E},d\mu_{M})$ the space of 
$L^{2}$ holomorphic sections of $E$ over $M$ with respect to $h_{E}$ and 
$d\mu_{M}$. 
Let $S$ be a closed  complex submanifold of $M$ and let $d\mu_{S}$ 
be a positive measure on $S$. 
The measured submanifold $(S,d\mu_{S})$ is said to be a set of 
interpolation for $(E,h_{E},d\mu_{M})$, or for the 
sapce $A^{2}(M,E,h_{E},d\mu_{M})$, if there exists a bounded linear operator
\[
I : A^{2}(S,E\mid_{S},h_{E},d\mu_{S})\longrightarrow A^{2}(M,E,h_{E},d\mu_{M})
\]
such that $I(f)\mid_{S} = f$ for any $f$. 
$I$ is called an interpolation operator.
The following theorem is crucial.

\begin{theorem}(\cite[Theorem 4]{o})
Let $M$ be a complex manifold with a continuous volume form $dV_{M}$,
let $E$ be a holomorphic vector bundle over $M$ with $C^{\infty}$-fiber 
metric $h_{E}$, let $S$ be a closed complex submanifold of $M$,
let $\Psi\in \sharp (S)$ and let $K_{M}$ be the canonical bundle of $M$.
Then $(S,dV_{M}(\Psi ))$ is a set of interpolation 
for $(E\otimes K_{M},h_{E}\otimes (dV_{M})^{-1},dV_{M})$, if 
the followings are satisfied.
\begin{enumerate}
\item There exists a closed set $X\subset M$ such that 
\begin{enumerate}
\item $X$ is locally negligble with respect to $L^{2}$-holomorphic functions, i.e., 
for any local coordinate neighbourhood $U\subset M$ and for any $L^{2}$-holomorphic function $f$ on $U\backslash X$, there exists a holomorphic function 
$\tilde{f}$ on $U$ such that $\tilde{f}\mid U\backslash X = f$.
\item $M\backslash X$ is a Stein manifold which intersects with every component of $S$. 
\end{enumerate}
\item $\Theta_{h_{E}}\geqq 0$ in the sense of Nakano,
\item $\Psi \in \sharp (S)\cap C^{\infty}(M\backslash S)$,
\item $e^{-(1+\epsilon )\Psi}\cdot h_{E}$ has semipositive 
curvature in the sense of Nakano for every $\epsilon \in [0,\delta]$ 
for some $\delta > 0$.
\end{enumerate}
Under these conditions, there exists a constant $C$ and an interpolation operator 
from $A^{2}(S,E\otimes K_{M}\mid_{S},h\otimes (dV_{M})^{-1}\mid_{S},dV_{M}[\Psi ])$
to $A^{2}(M,E\otimes K_{M},h\otimes (dV_{M})^{-1}.dV_{M})$ whose 
norm does not exceed $C\delta^{-3/2}$.
If $\Psi$ is plurisubharmonic, the interpolation operator can be chosen 
so that its norm is less than $2^{4}\pi^{1/2}$.
\end{theorem}

The above theorem can be generalized to the case that 
$(E,h_{E})$ is a singular hermitian line bundle with semipositive
curvature current  (we call such a singular hermitian line 
bundle $(E,h_{E})$ a {\bf pseudoeffective singular hermitian line bundle}) as was remarked in \cite{o}. 

\begin{lemma}
Let $M,S,\Psi ,dV_{M}, dV_{M}[\Psi], (E,h_{E})$ be as in Theorem 4.2. 
Let $(L,h_{L})$ be a pseudoeffective singular hermitian line 
bundle on $M$. 
Then $S$ is a set of interpolation for 
$(K_{M}\otimes E\otimes L,dV_{M}^{-1}\otimes h_{E}\otimes h_{L})$.  
\end{lemma}

\subsection{A construction of the function $\Psi$}
Let $M$ be a smooth projective $n$-fold and 
let $S$ be a $k$-dimensional (not necessary smooth)
subvariety of $M$. 
Let ${\cal U} = \{ U_{\gamma}\}$ be a finite 
Stein covering of $M$ and 
let $\{ f^{(\gamma )}_{1},\ldots , f_{m(\gamma )}^{(\gamma )}\}$ 
be a gnerator of the ideal sheaf associated with $S$ 
on $U_{\gamma}$. 
Let $\{ \phi_{\gamma}\}$ be a partition of unity subordinates
to ${\cal U}$. 
We set 
\[
\Psi := (n-k)\sum_{\gamma}\phi_{\gamma}\cdot (\sum_{\ell = 1}^{m(\gamma )}
\mid f_{\ell}^{(\gamma )}\mid^{2}).
\]
Then the residue volume form $dV[\Psi ]$ is defined 
as in the last subsection. 
Here the residue volume form $dV[\Psi ]$  
of an continuous volume form $dV$ on $M$ is not well defined on the singular
locus of $S$.
But this is not a difficulty to apply Theorem 4.3 or 
Lemma 4.1, since there exists a proper Zariski closed subset  
$Y$ of $X$ such that 
$(X - Y)\cap S$ is smooth. 

\subsection{Volume of pseudoeffective line bundles}

To measure the positivity of big line bundles on a projective 
variety, we shall introduce the notion of volume of a projective 
variety with respect to a big line bundle. 

\begin{definition} Let $L$ be a line bundle on a compact complex 
manifold $M$ of dimension $n$. 
We define the $L$-volume of $M$ by
\[
\mu (M,L) := n!\cdot\overline{\lim}_{m\rightarrow\infty}m^{-n}
\dim H^{0}(M,{\cal O}_{M}(mL)).
\]
\end{definition}
With respect to a pseudoeffective singular hermitian line bundle
(cf. for the definition of  pseudoeffective singular hermitian line bundles the last part of 4.2), 
we define the volume as follows. 
\begin{definition}(\cite{tu3})
Let $(L,h)$ be a pseudoeffective singular hermitian line bundle on a smooth projective variety
$X$ of dimension $n$. 
We define the volume $\mu (X,L)$ of $X$ with respect to 
$(L,h)$ by 
\[
\mu (X,(L,h)) := n!\cdot\overline{\lim}_{m\rightarrow\infty}m^{-n}
\dim H^{0}(X,{\cal O}_{X}(mL)\otimes{\cal I}(h^{m})).
\]
A pseudoeffective singular hermitian line bundle $(L,h)$ 
is said to be big, if $\mu (X,(L,h)) > 0$ holds. 

We may consider $\mu (X,(L,h))$ as the intersection number 
$(L,h)^{n}$. 
Let  $Y$ be  a subvariety of $X$ of dimension $d$
and let $\pi_{Y} : \tilde{Y}\longrightarrow Y$ be 
a resolution of $Y$.   We define 
$\mu (Y,(L,h)\mid_{Y})$ as 
\[
\mu (Y,(L,h)\mid_{Y})
:= \mu (\tilde{Y},\pi^{*}_{Y}(L,h)).
\]
The righthandside is independent of the choice of the 
resolution $\pi$ because of the remark below. 
\end{definition}
\begin{remark}
In Definition 4.3, let $\pi : \tilde{X}\longrightarrow X$ 
be any modification. 
Then 
\[
\mu (X,(L,h)) = \mu (\tilde{X},\pi^{*}(L,h))
\]
holds, since
\[
\pi_{*}({\cal O}_{\tilde{X}}(K_{\tilde{X}})\otimes {\cal I}(\pi^{*}h^{m}))
= {\cal O}_{X}(K_{X})\otimes {\cal I}(h^{m})
\]
holds for every $m$ and 
\[
\overline{\lim}_{m\rightarrow\infty}m^{-n}
\dim H^{0}(X,{\cal O}_{X}(mL)\otimes{\cal I}(h^{m}))
=  
\overline{\lim}_{m\rightarrow\infty}m^{-n}
\dim H^{0}(X,{\cal O}_{X}(mL + D)\otimes{\cal I}(h^{m}))
\]
holds for any Cartier divisor $D$ on $X$. 
This last equality can be easily checked, if $D$ is a smooth 
irreducible divisor, by using 
the exact sequence 
\[
0 \rightarrow {\cal O}_{X}(mL)\otimes {\cal I}(h^{m})
\rightarrow {\cal O}_{X}(mL+D)\otimes {\cal I}(h^{m})
\rightarrow {\cal O}_{D}(mL+D)\otimes {\cal I}(h^{m})
\rightarrow 0.
\]
For a general $D$, the equality follows by expressing $D$ 
as a difference of two very ample divisors. 
\end{remark}
\subsection{Construction of stratifications}

Let $X$ be a smooth projective $n$-fold of general type.
Let $h$ be an AZD of $K_{X}$ constructed as in Section 4.1. 
We may assume that $h$ is lowersemicontinuous (cf. \cite{tu2,d-p-s}). 
This is a technical assumption so that a local potential 
of the curvature current of $h$ is plurisubharmonic. 
This is used to restrict $h$ to a subvariety of $X$ 
(if we only assume that the  local potential is only locally integrable, 
the restriction is not well defined). 
We set 
\[
X^{\circ} = \{ x\in X\mid x\not{\in} \mbox{Bs}\mid mK_{X}\mid \mbox{and  
$\Phi_{\mid mK_{X}\mid}$ is a biholomorphism} 
\]
\[
\hspace{50mm} \mbox{on a neighbourhood of $x$ for some $m \geqq 1$}\} 
\]
as before. 
We set 
\[
\mu_{0} := \mu (X,(K_{X},h)) = \mu (X,K_{X}).
\]
The last equality holds, since $h$ is an AZD of $K_{X}$.
We note that for every $x\in X^{\circ}$, 
${\cal I}(h^{m})_{x}\simeq {\cal O}_{X,x}$ holds for every $m\geqq 0$ 
(cf. \cite{tu2} or \cite[Theorem 1.5]{d-p-s}). 
Using this fact the proof of the following lemma is identical
to that of Lemma 4.1. 

\begin{lemma} Let $x,x^{\prime}$ be distinct points on $X^{\circ}$.  
We set 
\[
{\cal M}_{x,x^{\prime}} = {\cal M}_{x}\otimes{\cal M}_{x^{\prime}},
\]
 where ${\cal M}_{x},{\cal M}_{x^{\prime}}$ denote the
maximal ideal sheaf of the points $x,x^{\prime}$ respectively.
Let $\varepsilon$ be a sufficiently small positive number.
Then 
\[
H^{0}(X,{\cal O}_{X}(mK_{X})\otimes{\cal M}_{x,x^{\prime}}^{\lceil\sqrt[n]{\mu_{0}}
(1-\varepsilon )\frac{m}{\sqrt[n]{2}}\rceil})\neq 0
\]
for every sufficiently large $m$.
\end{lemma}
Let us take a sufficiently large positive integer $m_{0}$ and let $\sigma_{0}$
be a general (nonzero) element of  
$H^{0}(X,{\cal O}_{X}(m_{0}K_{X})\otimes
{\cal M}_{x,x^{\prime}}^{\lceil\sqrt[n]{\mu_{0}}(1-\varepsilon )\frac{m_{0}}{\sqrt[n]{2}}\rceil})$.
We set 
\[
D_{0} := \frac{1}{m_{0}}(\sigma_{0})
\]
and  
\[
h_{0} = \frac{1}{\mid \sigma_{0}\mid^{2/m_{0}}}. 
\]
Let us define a positive (rational) number 
$\alpha_{0}$ by 
\[
\alpha_{0} := \inf\{\alpha > 0\mid 
\mbox{$\mbox{Supp}({\cal O}_{X}/{\cal I}(h_{0}^{\alpha}))$ contains  
both $x$ and $x^{\prime}$}\} .
\]
This is  essentially the same definition as in the last section. 
Then as before we see that 
\[
\alpha_{0}\leqq \frac{n\sqrt[n]{2}}{\sqrt[n]{\mu_{0}}(1-\varepsilon )}
\]
holds.
Suppose that for every small positive number $\delta$, 
$\mbox{Supp}({\cal O}_{X}/{\cal I}(h_{0}^{\alpha -\delta}))$ 
does not contain both $x$ and $x^{\prime}$. 
Other cases will be treated as before. 
Let $X_{1}$ be the minimal center of log canonical singularities
at $x$. 
We define the positive number $\mu_{1}$ by 
\[
\mu_{1}:= \mu (X_{1},(K_{X},h)\mid_{X_{1}}).
\]
Then since $x,x^{\prime}\in X^{\circ}$, 
$\mu_{1}$ is positive. 

For the later purpose, we shall perturb $h_{0}$ 
so that $X_{1}$ is the only center of log canonical singularities 
at $x$. 
Let $E$ be an effective {\bf Q}-divisor such that 
$K_{X} - E$ is ample. 
By the definition of $X^{\circ}$, we may assume that the support of $E$ does not 
contain $x$. 
Let $a$ be a positive integer such that 
$A:= a (K_{X}- E)$ is a very ample Cartier divisor 
such that ${\cal O}_{X}(A)\otimes {\cal I}_{X_{1}}$ is globally 
generated. 
Let $\rho_{1}, \ldots ,\rho_{e} 
\in H^{0}(X,{\cal O}_{X}(A)\otimes {\cal I}_{X_{1}})$ be a set of generators 
of ${\cal O}_{X}(A)\otimes {\cal I}_{X_{1}}$ on $X$. 
Then if we replace $h_{0}$ by 
\[
\frac{1}{(\mid \sigma_{0}\mid^{2}(\sum_{i=1}^{e}\mid\rho_{i}\mid^{2}))^{\frac{1}{m_{0}+a}}}
\]
has the desired property. 
If we take $m_{0}$ very large (in comparison with $a$), we can make 
the new $\alpha_{0}$ arbitrary close to the original $\alpha_{0}$.
To proceed further we use essentially the same procedure 
as in the previous section. 
The only essential difference here is to use Lemma 4.1 instead of Lemma 3.3. 
Let $x_{1},x_{2}$ be two distinct regular points on $X_{1}\cap X^{\circ}$. 
Let $m_{1}$ be a sufficiently large positive integer and 
Let 
\[
\sigma_{1}^{\prime}
\in 
H^{0}(X_{1},{\cal O}_{X_{1}}(m_{1}K_{X})\otimes {\cal I}(h^{m_{1}})
\cdot{\cal M}_{x_{1},x_{2}}^{\lceil\sqrt[n_{1}]{\mu_{1}}(1-\varepsilon )\frac{m}{\sqrt[n_{1}]{2}}\rceil})
\]
be a nonzero element. 

By Kodaira's lemma there is an effective {\bf Q}-divisor $E$ such
that $K_{X}- E$ is ample.
By the definition of $X^{\circ}$, we may assume that 
the support of $E$ does not contain both $x$ and $x^{\prime}$. 
Let $\ell_{1}$ be a sufficiently large positive integer which will be specified later  such that
\[
L_{1} := \ell_{1}(K_{X}- E)
\]
is Cartier.  Let $h_{L_{1}}$ be a $C^{\infty}$-hermitian metric on $L_{1}$ 
with strictly positive curvature. 
Let $\tau$ be a nonzero section in 
$H^{0}(X,{\cal O}_{X}(L_{1}))$.  
We set 
\[
\Psi := \alpha_{0} \log\frac{h_{0}}{h}.
\]
Let $dV$ be a $C^{\infty}$-volume form on $X$.
We note that the residue volume form  $dV [\Psi ]$ on $X_{1}$
 may have pole along some proper subvarieties in $X_{1}$. 
By taking $\ell_{1}$ sufficiently large and taking $\tau$ properly,
we may assume that $h_{L_{1}}(\tau ,\tau )\cdot dV[\Psi ]$ is nonsingular on 
$X_{1}$ in the sense that the pullback of  it to a nonsingular model of $X_{1}$
is nonsingular. 
Then by applying  Lemma 4.1 for $(E,h_{E})$ to be 
$((\lceil 1+\alpha_{0}\rceil )K_{X},h^{\lceil 1+\alpha_{0}\rceil})$ 
and 
$(X,X_{1},\Psi ,dV,dV[\Psi], ((m_{1}-\lceil \alpha_{0}\rceil -2)K_{X}+L_{1},h^{m_{1}-1}h_{L_{1}})$ 
we see that   
\[
\sigma_{1}^{\prime}\otimes\tau\in
H^{0}(X_{1},{\cal O}_{X_{1}}(m_{1}K_{X}+L_{1})
\otimes {\cal I}(h^{m_{1}})\cdot{\cal M}_{x_{1},x_{2}}^{\lceil\sqrt[n_{1}]{\mu_{1}}(1-\varepsilon )\frac{m_{1}}{\sqrt[n_{1}]{2}}\rceil})
\]
extends to a section
\[
\sigma_{1}\in H^{0}(X,{\cal O}_{X}((m_{1}+\ell_{1} )K_{X})).
\]
We note that even though $dV[\Psi ]$ may have singularity on $X_{1}$, 
we may apply
Lemma 4.1, because there exists a proper Zariski closed subset  
$Y$ of $X$ such that the restriction of $dV[\Psi ]$
to $(X - Y)\cap X_{1}$ is smooth. 
Of course the singularity of $dV[\Psi ]$ affects to the $L^{2}$-condition.
But this has already been handled by the boundedness of 
$h_{L_{1}}(\tau ,\tau )\cdot dV[\Psi ]$. 

Then by entirely the same procedure as in Section 3,
 for distinct points $x,x^{\prime}$, we construct 
a strictly decreasing sequence of subvarieties 
\[
X = X_{0}\supset X_{1}\supset \cdots \supset X_{r}\supset X_{r+1} =  x\,\,\mbox{or}\,\, x^{\prime}
\]
and invariants (depending on small positive numbers $\varepsilon_{0},\ldots ,
\varepsilon_{r-1}$, large positive integers $m_{0},m_{1},\ldots ,m_{r}$, etc.) :
\[
\alpha_{0} ,\alpha_{1},\ldots ,\alpha_{r},
\]
\[
\varepsilon_{0},\varepsilon_{1},\ldots ,\varepsilon_{r-1}
\]
\[
\mu_{0},\mu_{1},\ldots ,\mu_{r}
\]
and
\[
n >  n_{1}> \cdots > n_{r}
\]
inductively.
By Nadel's vanishing theorem (\cite[p.561]{n}) we have the following lemma.

\begin{lemma} 
Let $x,x^{\prime}$ be two distinct points on $X^{\circ}$. 
Then for every $m\geqq \lceil\sum_{i=0}^{r}\alpha_{i}\rceil +1$,
$\Phi_{\mid mK_{X}\mid}$ separates $x$ and $x^{\prime}$.
\[
\alpha_{i}\leqq \frac{n_{i}\sqrt[n_{i}]{2}}{\sqrt[n_{i}]{\mu_{i}}} + \delta 
\]
hold for $1\leqq i\leqq r$.
\end{lemma}

We may also construct the stratification as a family as in 
Section 3.2. 

\subsection{Another subadjunction theorem}
Let $M$ be a smooth projective variety 
and let $(L,h_{L})$ be a singular hermitian line bundle on $M$ such that 
$\Theta_{h_{L}}\geqq 0$ on $M$.  
Let $dV$ be a $C^{\infty}$-volume form on $M$. 
Let $\sigma \in \Gamma (\bar{M},{\cal O}_{\bar{M}}(m_{0}L)\otimes {\cal I}(h))$ be a 
global section. 
Let $\alpha$ be a positive rational number $\leqq 1$ and let $S$ be 
an irreducible subvariety of $M$ 
such that  $(M, \alpha (\sigma ))$ is log canonical but not KLT(Kawamata log-terminal)
on the generic point of $S$ and $(M,(\alpha -\epsilon )(\sigma ))$ is KLT on the generic point of $S$ 
for every $0 < \epsilon << 1$. 
We set 
\[
\Psi = \alpha \log h_{L}(\sigma ,\sigma ).
\]
Suppose that $S$ is smooth 
(if $S$ is not smooth, we just need to take an embedded 
resolution and apply Theorem 4.4 below). 
We shall assume that $S$ is not contained in the 
singular locus of $h$, where the singular locus of $h$ means the 
set of points where $h$ is $+\infty$. 
Then as in Section 4.2, we may define a (possibly singular measure) 
$dV[\Psi ]$ on $S$. 
This can be viewed as follows. 
Let $f : N \longrightarrow M$ be a log-resolution of 
$(X,\alpha (\sigma ))$. 
Then  as before we may define the singular volume form $f^{*}dV[f^{*}\Psi ]$ 
on the divisorial component of $f^{-1}(S)$. 
The singular volume form $dV[\Psi ]$ is defined as the fibre integral of 
$f^{*}dV[f^{*}\Psi ]$. 
Let $d\mu_{S}$ be a $C^{\infty}$-volume form on $S$ and 
let $\varphi$ be the function on $S$ defined by
\[
\varphi := \log \frac{dV[\Psi ]}{d\mu_{S}}
\]
($dV[\Psi ]$ may be singular on a subvariety of $S$, also 
it may be totally singular on $S$). 

\begin{theorem}(\cite[Theorem 5.1]{tu5})
Let $M$,$S$,$\Psi$ be as above. 
Suppose that $S$ is smooth.   
Let $d$ be a positive integer such that $d \geqq \alpha m_{0}$. 
Then every element of 
$A^{2}(S,{\cal O}_{S}(m(K_{M}+dL)),e^{-(m-1)\varphi}\cdot dV^{-m}\otimes 
h_{L}^{m}\mid_{S},dV[\Psi ])$ 
extends to an element of 
\[
H^{0}(M,{\cal O}_{M}(m(K_{M}+dL))). 
\]
\end{theorem}
As we mentioned as above the smoothness assumption on $S$ is 
just to make the statement simpler.  

Theorem 4.5 follows from 
Theorem 4.6 below by  minor modifications (cf. \cite{tu5}). 
The main difference is the fact that the residue volume form 
$dV[\Psi ]$ is singular on $S$. 
But this does not affect the proof, since in the $L^{2}$-extension theorem
(Theorem 4.4) we do not need to assume that the manifold $M$ is compact. 
Hence we may remove a suitable subvarieties so that we do not need 
to consider the pole of $dV[\Psi]$ on $S$ (but of course the pole of 
$dV[\Psi]$  affects the $L^{2}$-conditions). 

\begin{theorem}
Let $M$ be a projective manifold with a continuous volume form $dV$,
let $L$ be a holomorphic line bundle over $M$ with a $C^{\infty}$-hermitian metric  $h_{L}$, let $S$ be a compact complex submanifold of $M$,
let $\Psi : M \longrightarrow [-\infty ,0)$ be a continuous function and let $K_{M}$ be the canonical bundle of $M$.
\begin{enumerate}
\item $\Psi \in  \sharp (S) \cap C^{\infty}(M\backslash S)$ (As for the 
definition of $\sharp (S)$, see Section 4.2), 
\item $\Theta_{h\cdot e^{-(1+\epsilon )\Psi}}\geqq 0$ for 
every $\epsilon \in [0,\delta ]$ for some $\delta > 0$,
\item there is a positive line bundle on $M$.
\end{enumerate}
Then every element of  $H^{0}(S,{\cal O}_{S}(m(K_{M}+L)))$ extends to an element of 
$H^{0}(M,{\cal O}_{M}(m(K_{M}+L)))$. 
\end{theorem}
For the completeness we shall give a proof of Theorem 4.5 
under the additional condition :  \vspace{5mm} \\
{\bf Condition} $(K_{M}+L,dV^{-1}\otimes h_{L})$ 
is big (cf. Definition 4.3). 
\vspace{5mm} \\
The reason why we put this condition is that we only need Theorem 4.5 and 4.6 under this condition. 
 
Let $M, S, L$ be as in Theorem 4.6.
Let $n$ be the dimension of $M$.  Let $h_{S}$ be a canonical AZD
(\cite{tu2}) of 
$K_{M}+L\mid_{S}$. 
Let $A$ be a sufficiently ample line bundle on $M$. 
We consider the Bergman kernel 
\[
K(S,A+m(K_{M}+L)\mid_{S},h_{A}\cdot h_{S}^{m-1}\cdot dV^{-1}\cdot 
h_{L},d\Psi_{S}) = \sum_{i}\mid \sigma_{i}^{(m)}\mid^{2}, 
\]
where $\{ \sigma_{i}^{(m)}\}$ is a complete orthonormal basis of 
$A^{2}(S,A+m(K_{M}+L)\mid_{S},h_{A}\cdot h_{S}^{m-1}\cdot dV^{-1}\cdot 
h_{L},d\Psi_{S})$. 
We note that (cf. \cite[p.46, Proposition 1.4.16]{kr})
\[
K(S,A+m(K_{M}+L)\mid_{S},h_{A}\cdot h_{S}^{m-1}\cdot dV^{-1}\cdot 
h_{L},d\Psi_{S})(x)
\]
\[= \sup \{ \mid\sigma\mid^{2}(x) \, \mid\, 
\sigma \in A^{2}(S,A+m(K_{M}+L)\mid_{S},h_{A}\cdot h_{S}^{m-1}\cdot dV^{-1}\cdot 
h_{L},d\Psi_{S}), \parallel \sigma\parallel = 1\}  
\]
holds for every $x\in S$.
 
Let us define the singular hemitian metric on $m(K_{M}+L)\mid_{S}$ by 
\[
h_{m,S} := K(A+m(K_{M}+L)\mid_{S},h_{A}\cdot h_{S}^{m-1}\cdot dV^{-1}\cdot 
h_{L},d\Psi_{S})^{-1}
\]
Then as in \cite[Section 4]{d} and \cite{tu2}, we see that 
\[
h_{S} : = \liminf_{m\rightarrow\infty}\sqrt[m]{h_{m,S}}
\]
holds.   Hence $\{\sqrt[m]{h_{m,S}}\}$  is considered to be an 
algebraic approximation of $h_{S}$.
We note that there exists a positive constant $C_{0}$ independent of $m$ such that 
\[
h_{m,S} \leqq C_{0}\cdot h_{A}\cdot h_{S}^{m}
\]
holds for every $m\geqq 1$ as in \cite{tu2}.
Let $h_{M}$ be a  canonical AZD of $K_{M}+L$.
By the assumption, $(K_{M}+L,h_{M})$ is big , i.e., 
\[
 \lim_{m\rightarrow\infty}
\frac{\log \dim H^{0}(M,{\cal O}_{M}(A+m(K_{M}+L))
\otimes {\cal I}(h_{M}^{m}))}{\log m} = n.  
\]
Let $\nu(S)$ denote the numerical Kodaira dimension of 
$(K_{M}+L\mid_{S},h_{S})$ (for our purpose we may assume that 
$\nu (S) = \dim S$), i.e.,
\[ 
\nu (S):= \lim_{m\rightarrow\infty}
\frac{\log \dim H^{0}(S,{\cal O}_{S}(A+m(K_{M}+L))
\otimes {\cal I}(h_{S}^{m}))}{\log m}.  
\]
Inductively on $m$, we extend each 
\[
\sigma\in
A^{2}(S,A+m(K_{M}+L)\mid_{S},h_{A}\cdot h_{S}^{m-1}\cdot dV^{-1}\cdot 
h_{L},d\Psi_{S}) 
\]
to a section 
\[
\tilde{\sigma}\in A^{2}(M,A+ m(K_{M}+L),dV^{-1}\cdot h_{L}\cdot \tilde{h}_{m-1},dV)
\]
with the estimate
\[
\parallel \tilde{\sigma}\parallel\leqq C\cdot m^{-(n-\nu(S))} \parallel \sigma\parallel ,
\]
where  
$\parallel\,\,\,\,\parallel$'s denote the $L^{2}$-norms respectively, 
$C$ is a positive constant indpendent of $m$ 
and we have defined 
\[
\tilde{K}_{m}(x) := \sup\{ \mid \tilde{\sigma}\mid^{2}(x) \mid \,\,
\parallel\tilde{\sigma}\mid_{S}\parallel = 1, 
\parallel\tilde{\sigma}\parallel \leqq C\cdot m^{-(n-\nu (S))} \}
\]
and set
\[
\tilde{h}_{m} = \frac{1}{\tilde{K}_{m}}.  
\]
If we take $C$ sufficiently large, then $\tilde{h}_{m}$ 
is well defined for every $m\geqq 0$. 
Here we note that the factor $m^{-(n-\nu (S))}$ comes from that 
$h_{M}$ is dominated by a singular hermitian metric with 
strictly positive curvature by Kodaira's lemma by the condition 
(cf.  \cite[p,105,(1,11)]{ti}).
By easy inductive estimates, we see that 
\[
\tilde{h}_{\infty} :=  \liminf_{m\rightarrow\infty}\sqrt[m]{\tilde{h}_{m}}
\]
exists and gives an extension of $h_{S}$.
In fact $h_{A}^{1/m}(\tilde{K}_{m})^{1/m}$ is 
uniformly bounded from above as in \cite[p.127,Lemma 3.3]{tu6}.
Then by \cite{o} for every $m\geqq 1$, we may extend every element of 
$A^{2}(m(K_{M}+L)\mid_{S},dV^{-1}\cdot h_{L}\cdot h_{S}^{m-1},dV[\Psi ])$ to 
$A^{2}(m(K_{M}+L),dV^{-1}\cdot h_{L}\cdot \tilde{h}_{\infty},dV)$. 
This completes the proof of Theorem 4.6. 

\subsection{Positivity result}

The following positivity theorem is a key to the proof of Theorem 1.1 and 
1.2.
\begin{theorem}(\cite[Theorem 2]{ka})
Let $f : X \longrightarrow B$ be a surjective morphism of smooth projective 
varieties with connected fibers.
Let $P = \sum P_{j}$ and $Q = \sum_{\ell}Q_{\ell}$ be normal crossing divisors on $X$ and $B$ respectively, such that $f^{-1}(Q) \subset P$ and $f$ 
is smooth over $B\backslash Q$.
Let $D = \sum d_{j}P_{j}$ be a {\bf Q}-divisor on $X$, where $d_{j}$ may be positive, zero or negative, which satisfies the following conditions :
\begin{enumerate}
\item $D = D^{h} + D^{v}$ such that 
$f :\mbox{Supp}(D^{h})\rightarrow B$ is surjective and smooth over $B\backslash Q$, and $f(\mbox{Supp}(D^{h}) \subset Q$.
An irreducible component of $D^{h}$(resp. $D^{v}$) is called horizontal
(resp. vertical).
\item $d_{j} < 1$ for all $j$.
\item The natural homomorphism ${\cal O}_{B}\rightarrow f_{*}{\cal O}_{X}(\lceil -D\rceil )$ is surjective at the generic point of $B$.
\item $K_{X} +  D\sim_{\mbox{\bf Q}}f^{*}(K_{B} + L)$ for some 
{\bf Q}-divisor $L$ on $B$.
Let 
\begin{eqnarray*}
f^{*}Q_{\ell}& =  &\sum_{j}w_{\ell j}P_{j} \\
\bar{d}_{j} & :=  & \frac{d_{j} +w_{\ell j}-1}{w_{\ell j}}\,\,\,\,\mbox{if}\,\,\,\,
f(P_{j}) = Q_{\ell} \\
\delta_{\ell} &: =  & \max \{\bar{d}_{j} ; f(P_{j}) = Q_{\ell}\} \\
\Delta & :=  & \sum_{\ell}\delta_{\ell}Q_{\ell} \\
M & :=  & L - \Delta .
\end{eqnarray*}
Then $M$ is nef.
\end{enumerate} 
\end{theorem} 
 
Here the meaning of the divisor $\Delta$ may be difficult to understand.
I would like to give an geometric interpretation of $\Delta$.  
Let $X,P,Q,D,B,\Delta$ be as above. Let $dV$ be a 
$C^{\infty}$-volume form on $X$. 
Let $\sigma_{j}$ be a global section of ${\cal O}_{X}(P_{j})$
with divisor $P_{j}$. 
Let $\parallel\sigma_{j}\parallel$ denote the hermitian norm 
of $\sigma_{j}$ with respect to a $C^{\infty}$-hermitian metric
on ${\cal O}_{X}(P_{j})$ respectively. 
Let us consider the singular volume form
\[
\Omega := \frac{dV}{\prod_{j}\parallel\sigma_{j}\parallel^{2d_{j}}}.
\]
Then by taking the fiber integral of $\Omega$ with respect to 
$f : X \longrightarrow B$, we obtain a singular volume form 
$\int_{X/B}\Omega$ on $B$. 
Then the divisor $\Delta$ corresponds to the singularity 
of the singular volume form $\int_{X/B}\Omega$ on $B$. 

\subsection{Use of two subadjunction theorems}

This subsection is the counterpart of Section 3.3. 
\begin{lemma}
For every  $X_{j}$, 
\[
\mu (W_{j},K_{W_{j}})
\leqq 
(\lceil (1+\sum_{i=0}^{j-1}\alpha_{i})\rceil )^{n_{j}}\cdot\mu_{j}
\]
holds, where $W_{j}$ is a nonsingular model of $X_{j}$ (we note that
 $\mu (W_{j},K_{W_{j}})$ is independent of the choice of the nonsingular model 
 $W_{j}$).
\end{lemma}
{\bf Proof}.

Let us set 
\[
\beta_{j}:= \varepsilon_{j-1}+\sum_{i=0}^{j-1} (\alpha_{i}-\varepsilon_{i}).
\]
Let $D_{i}$ denotes the divisor $m_{i}^{-1}(\sigma_{i})$ and we set 
\[
D : = \sum_{i=1}^{j-1} (\alpha_{i}-\varepsilon_{i})D_{i}+\varepsilon_{j-1}D_{j-1}.
\]
Let $\pi : Y \longrightarrow X$ be a log resolution of 
$(X,D)$ 
which factors through the embedded resolution 
$\varpi : W_{j} \longrightarrow X_{j}$ of $X_{j}$. 
By the perturbation as in Section 4.5, we may assume that there exists a unique  irreducible component $F_{j}$ of the exceptional divisor with discrepancy $-1$ which dominates $X_{j}$. 
Let 
\[
\pi_{j} : F_{j} \longrightarrow W_{j}
\]
be the natural morphism induced by the construction. 
We set 
\[
\pi^{*}(K_{X} + D)\mid_{F_{j}} = K_{F_{j}} + G.
\]
We may assume that the support of $G$ is a divisor with normal crossings.
Then all the coefficients of the horizontal component $G^{h}$  with respect to $\pi_{j}$ are less than $1$  because $F_{j}$ is the unique exceptional divisor with discrepancy $-1$. 

Let $dV$ be a $C^{\infty}$-volume form on the  $X$. 
Let $\Psi$ be the function defined by 
\[
\Psi := \log (h^{\beta_{j}}\cdot\mid\sigma_{j-1}\mid^{{\frac{2\varepsilon_{j-1}}{m_{j-1}}}}\cdot\prod_{i=0}^{j-1}\mid \sigma_{i}\mid^{\frac{2(\alpha_{i}-\varepsilon_{i})}{m_{i}}}).
\]
Then the residue $\mbox{Res}_{F_{j}}(\pi^{*}(e^{-\Psi}\cdot dV))$ of  $\pi^{*}(e^{-\Psi}\cdot dV)$
to $F_{j}$ 
is a singular volume form with algebraic singularities corresponding 
to the divisor $G$. 
Since every coefficient of $G^{h}$ is less than $1$ , there exists a Zariski open subset 
$W_{j}^{0}$ of $W_{j}$ such that 
 $\mbox{Res}_{F_{j}}(\pi^{*}(e^{-\Psi}\cdot dV))$ is integrable
 on $\pi_{j}^{-1}(W_{j}^{0})$.

Then  the pullback of the residue $dV[\Psi ]$ of $e^{-\Psi}\cdot dV$ 
(to $X_{j}$) to 
 $W_{j}$ is given by the fiber integral of 
the above singular volume form $\mbox{Res}_{F_{j}}(\pi^{*}(e^{-\Psi}\cdot dV))$ on $F_{j}$, i.e.,
\[
\varpi^{*}dV[\Psi ] = \int_{F_{j}/W_{j}}\mbox{Res}_{F_{j}}(\pi^{*}(e^{-\Psi}\cdot dV))
\]
holds. 
By Theorem 4.7, we see that 
$(K_{F_{j}}+G) -\pi_{j}^{*}(K_{W_{j}}+\Delta )$ is nef,
where $\Delta$ is the divisor defined as in Theorem 4.7. 
And as in \cite[cf. Proof of Theorem 1]{ka}, $\varpi_{*}\Delta$ is effective 
on $X_{j}$. 
We may assume that $\pi : Y \longrightarrow X$ is a simultaneous resolution
of $\mbox{Bs}\mid m_{i}K_{X}\mid (0\leqq i\leqq j-1)$. 
Let us decompose $D_{i} (0\leqq i\leqq j-1)$ as 
\[
\pi^{*}D_{i} = P_{i} + N_{i},
\]
where $P_{i}$ is the free part and $N_{i}$ is the fixed component. 
Let us consider the contribution of 
\[
\sum_{i=0}^{j-1}\frac{\alpha_{i}-\varepsilon_{i}}{m_{i}}N_{i}
+ \frac{\varepsilon_{j-1}}{m_{j-1}}N_{j-1}
\]
to the divisor $\Delta$. 
Let $\sigma_{\Delta}$ be a multivalued meromorphic section of the 
{\bf Q}-line bundle ${\cal O}_{W_{j}}(\Delta )$ with divisor $\Delta$.
Let $h_{\Delta}$ be a $C^{\infty}$-hermitian metric on the {\bf Q}-line bundle
${\cal O}_{W_{j}}(\Delta )$. 
Then since
\[
\sigma_{i}\in H^{0}(X,{\cal O}_{X}(m_{i}K_{X})\otimes{\cal I}_{\infty}(h^{m_{i}}))
\]
for every $0\leqq i\leqq j-1$ (cf. Remark 4.3),
we see that 
\[
h_{\Delta}(\sigma_{\Delta},\sigma_{\Delta}) =  
O(\varpi^{*}((dV\cdot h^{-1})\mid_{X_{j}})^{\beta_{j}})
\]
holds. 
Hence we see that 
\begin{equation}
\mu (W_{j},K_{W_{j}}) \leqq \mu (X_{j},((\lceil 1+\beta_{j}\rceil )K_{X},
dV^{-1}\otimes h^{\lceil\beta_{j}\rceil})\mid_{X_{j}})
\end{equation}
holds. 

Let us recall the interpretation of the divisor $\Delta$ in 
Section 4.7.
Let $dV_{W_{j}}$ be a $C^{\infty}$-volume form on $W_{j}$. 
Then there exists a positive constant $C > 1$ such that 
\[
\varpi^{*}dV[\Psi ] = \int_{F_{j}/W_{j}}\mbox{Res}_{F_{j}}(\pi^{*}(e^{-\Psi}\cdot dV))
\leqq C\cdot\frac{\varpi^{*}(dV\cdot h^{-1})^{\beta_{j}}}{h_{\Delta}(\sigma_{\Delta},\sigma_{\Delta})}\cdot dV_{W}
\] 
hold.
We note that $dV\cdot h^{-1}$ is bounded from above by the construction of $h$,

Since  $\varpi_{*}\Delta$ is effective 
on $X_{j}$,  by applying  Theorem 4.5 , we have the interpolation :
\[
A^{2}(W_{j},m(\lceil 1+\beta_{j}\rceil )K_{X}, 
e^{-(m-1)\varphi}\cdot dV^{-m}\otimes h^{m\lceil \beta_{j}\rceil},
dV[\Psi ])
\rightarrow 
H^{0}(X,{\cal O}_{X}(m(\lceil 1+ \beta_{j}\rceil )K_{X})),
\]
where $\varphi$ is the weight function defined as in Theorem 4.5.

Since $h$ is an AZD of $K_{X}$ constructed as in Section 4.1, we see that every element of 
$H^{0}(X,{\cal O}_{X}(m(\lceil 1+ \beta_{j}\rceil )K_{X}))$ 
is  bounded on $X$  with respect to $h^{m(\lceil 1+\beta_{j}\rceil )}$
(cf. Remark 4.3). 
In particular the restriction of an element of $H^{0}(X,{\cal O}_{X}(m(\lceil 1+ \beta_{j}\rceil )K_{X}))$ to $X_{j}$ is bounded with respect to 
$h^{\lceil 1+\beta_{j}\rceil}\mid_{X_{j}}$.
Hence by the existence of the above extension, we have that 
\begin{equation}
\mu (W_{j},K_{W_{j}})
\leqq 
\mu (X_{j},(\lceil 1+\beta_{j}\rceil )K_{X},
h^{\lceil 1+\beta_{j}\rceil})\mid_{X_{j}})
\end{equation}
holds. 
The difference between the inequalities (1) and (2) is that 
in (2) $dV^{-1}\otimes h^{\lceil\beta_{j}\rceil}\mid_{X_{j}}$ is 
replaced by  $h^{\lceil 1+\beta_{j}\rceil}\mid_{X_{j}}$.
This is the only point where Theorem 4.5 is used. 

By the trivial inequality 
\[
\beta_{j} \leqq \sum_{i=0}^{j}\alpha_{i}.
\]
we have that 
\[
\mu (W_{j},K_{W_{j}})\leqq (\lceil 1+\sum_{i=0}^{j-1}\alpha_{i}\rceil )^{n_{j}}(K_{X},h)^{n_{j}}\cdot X_{j}\]
holds by the definition of $(K_{X},h)^{n_{j}}\cdot X_{j}$.
This is the desired inequality, since $\mu_{j} = (K_{X},h)^{n_{j}}\cdot X_{j}$ 
holds by definition.
{\bf Q.E.D.}

\subsection{Estimate of the degree}
\begin{lemma}
If $\Phi_{\mid mK_{X}\mid}$ is a birational rational map
onto its image, then
\[
\deg \Phi_{\mid mK_{X}\mid}(X)\leqq \mu_{0}\cdot m^{n}
\]
holds.
\end{lemma}
{\bf Proof}.
Let $p : \tilde{X}\longrightarrow X$ be the resolution of 
the base locus of $\mid mK_{X}\mid$ and let 
\[
p^{*}\mid mK_{X}\mid = \mid P_{m}\mid + F_{m}
\]
be the decomposition into the free part $\mid P_{m}\mid$ 
and the fixed component $F_{m}$. 
We have
\[
\deg \Phi_{\mid mK_{X}\mid}(X) = P_{m}^{n}
\]
holds.
Then by the ring structure of $R(X,K_{X})$, we have an injection 
\[
H^{0}(\tilde{X},{\cal O}_{\tilde{X}}(\nu P_{m}))\rightarrow 
H^{0}(X,{\cal O}_{X}(m\nu K_{X})\otimes{\cal I}(h^{m\nu}))
\]
for every $\nu\geqq 1$, since 
the righthandside is isomorphic to 
$H^{0}(X,{\cal O}_{X}(m\nu K_{X}))$ by the definition of 
an AZD.
We note that since ${\cal O}_{\tilde{X}}(\nu P_{m})$ is globally generated
on $\tilde{X}$, for every $\nu \geqq 1$ we have the injection 
\[
{\cal O}_{\tilde{X}}(\nu P_{m})\rightarrow p^{*}({\cal O}_{X}(m\nu K_{X})\otimes{\cal I}(h^{m\nu})).
\]
Hence there exists a natural homomorphism 
\[
H^{0}(\tilde{X},{\cal O}_{\tilde{X}}(\nu P_{m}))
\rightarrow 
H^{0}(X,{\cal O}_{X}(m\nu K_{X})\otimes{\cal I}(h^{m\nu}))
\]
for every $\nu\geqq 1$. 
This homomorphism is clearly injective. 
This implies that 
\[
\mu_{0} \geqq  m^{-n}\mu (\tilde{X},P_{m})
\]
holds by the definition of $\mu_{0}$. 
Since $P_{m}$ is nef and big on $X$, we see that 
\[
\mu (\tilde{X},P_{m}) = P_{m}^{n}
\]
holds.
Hence
\[
\mu_{0}\geqq m^{-n}P_{m}^{n}
\]
holds.  This implies that
\[
\deg \Phi_{\mid mK_{X}\mid}(X)\leqq \mu_{0}\cdot m^{n}
\]
holds.
{\bf Q.E.D.}
\subsection{Completion of the proof of Theorem 1.1 and 1.2.}

We shall complete the proofs of Theorem 1.1 and 1.2. 
Suppose that Theorem 1.2 holds for every projective varieties of general 
type of dimension $< n$, i.e., there exists a positive constants $\{C(k) (k <n)\}$
such that for every every smooth projective $k$-fold $Y$ of general type
\[
\mu (Y,K_{Y}) \geqq C(k)
\] 
holds. 
Let $X$ be a smooth projective variety of general type of 
dimension $n$. 
Let $U_{0}$ be a nonempty Zariski open subset of $X$ with respect to 
{\bf coutable Zariski topology} such that for every $x\in U_{0}$ 
there exist no subvarieties of nongeneral type containing $x$. 
Such $U_{0}$ surely exists, since there exists no dominant family of 
 subvarieties of nongeneral type in $X$ (see Section 3.3).  
Then if $(x,x^{\prime})\in U_{0}\times U_{0} -\Delta_{X}$, the stratum $X_{j}$ 
as in Section 3.1 is of general type as before. 
By Lemma 4.4, we see that
\[
C(n_{j}) \leqq (\lceil (1+\sum_{i=0}^{j-1}\alpha_{i})\rceil )^{n_{j}}\cdot \mu_{j}
\]
holds for $W_{j}$. 
As in Lemma 3.10 in  Section 3.3,we obtain the following lemma. 
\begin{lemma}
Suppose that $\mu_{0} \leqq 1$ holds.
Then there exists a positive constant $C$ depending only on $n$ 
such that for every $(x,x^{\prime})\in U_{0}\times U_{0} -\Delta_{X}$ 
the corresponding invariants $\{ \mu_{0},\cdots ,\mu_{r}\}$ 
and $\{ n_{1},\cdots ,n_{r}\}$ depending 
on $(x,x^{\prime})$ ($r$ may also depend on $(x,x^{\prime})$) satisfies
the inequality :
\[
1+\lceil \sum_{i=0}^{r}\frac{\sqrt[n_{i}]{2}\, n_{i}}{\sqrt[n_{i}]{\mu_{i}}}\rceil \leqq \lfloor\frac{C}{\sqrt[n]{\mu_{0}}}\rfloor .
\]
\end{lemma} 
We note that the finite stratification of $U_{0}\times U_{0} -\Delta_{X}$ 
is unnecessary for the proof, because $\{ n_{1},\cdots ,n_{r}\}$ is a 
strictly decreasing sequence and this sequence has only finitely many 
possibilities. 
By Lemma 4.3  we see that 
for 
\[
m := \lfloor\frac{C}{\sqrt[n]{\mu_{0}}}\rfloor ,
\]
$\mid mK_{X}\mid$ separates points in $U_{0}$.
Hence  $\mid mK_{X}\mid$ gives a birational embedding of $X$.

Then by Lemma 4.5, if $\mu_{0} \leqq 1$ holds, 
\[
\deg \Phi_{\mid mK_{X}\mid}(X) 
\leqq C^{n}
\]
holds. 
Also 
\[
\dim H^{0}(X,{\cal O}_{X}(mK_{X}))
\leqq n+1 + \deg \Phi_{\mid mK_{X}\mid}(X) 
\]
holds by the semipositivity of the $\Delta$-genus (\cite{fu}). 
Hence we have that if $\mu_{0} \leqq 1$,
\[
\dim H^{0}(X,{\cal O}_{X}(mK_{X})) 
\leqq n+1+ C^{n} 
\]
holds. 

Since $C$ is a positive constant depending only on $n$, 
combining  the above two inequalities, by the argument as in Section 3.5, 
we have that there exists a positive constant $C(n)$ 
depending only on $n$ such that 
\[
\mu_{0}  \geqq C(n)
\]
holds. 

Then as in Section 3.5 we see that there exists 
a positive integer $\nu_{n}$ depending only on $n$ such that 
for every projective $n$-fold $X$ of general type, 
$\mid mK_{X}\mid$ gives a birational embedding into a 
projective space for every $m\geq \nu_{n}$. 
This completes the proof of Theorem 1.1.

\section{The Severi-Iitaka conjecture}

Let $X$ be a smooth projective variety.
We set
\[
Sev (X) := \{ (f,[Y])\mid f : X\longrightarrow Y\,\,\,\,\mbox{dominant
rational map and $Y$ is of general type}\} ,
\]
where $[Y]$ denotes the birational class of $Y$.
By Theorem 1.1 and \cite[p.117, Proposition 6.5]{m}
we obtain the following theorem.

\begin{theorem}
$Sev (X)$ is finite.
\end{theorem}
\begin{remark}
In the case of $\dim Y = 1$, Theorem 5.1 is known 
as Severi's theorem. 
In the case of $\dim Y = 2$, Theorem 5.1 has already been known 
by K. Maehara (\cite{m}).
In the case of $\dim Y = 3$, Theorem 5.1 has recently proved by
T. Bandman and G. Dethloff (\cite{b-d}).
\end{remark}

\section{Appendix}
\subsection{Volume of nef and big line bundles}
The following fact seems to be well known. 
But for the completeness, I would like to include the proof. 
\begin{proposition}
Let $M$ be a smooth projective $n$-fold and let $L$ be 
a nef and big line bundle on $M$. 
Then 
\[
n!\cdot \overline{\lim}_{m\rightarrow \infty}
m^{-n}\dim H^{0}(M,{\cal O}_{M}(mL)) 
= L^{n}
\]
holds. 
\end{proposition}
{\bf Proof of Proposition 5.1}. 
Since $L$ is big, there exists an effective {\bf Q}-divisor 
$F$ such that $L - F$ is ample. 
Let $a$ be a positive integer such that 
$A := a(L - F)$ is a very ample Cartier divisor and 
$A - K_{X}$ is ample. 
Then by  the Kodaira vanishing theorem, for every 
$q\geqq 1$, 
\[
H^{q}(M,{\cal O}_{M}(A + mL)) = 0
\]
holds for every $m\geqq 0$. 
By the Riemann-Roch theorem we have that 
\[
n!\cdot \overline{\lim}_{m\rightarrow \infty}
m^{-n}\dim H^{0}(M,{\cal O}_{M}(A+ mL)) 
= L^{n}
\]
holds. 
By the definition of $A$, we see that 
\[
n!\cdot \overline{\lim}_{m\rightarrow \infty}
m^{-n}\dim H^{0}(M,{\cal O}_{M}(mL)) 
= L^{n}
\]
holds. 
This completes the proof. {\bf Q.E.D.}

\subsection{A Serre type vanishing theorem}
\begin{lemma}
Let $X$ be a projective variety with only canonical singularities
(cf. \cite[p.56, Definition 2.34]{k-m}). 
Let $E$ be a vector bundle on $X$ and let 
$L$ be a nef line bundle on $X$.  
Let $A$ be an ample line bundle on $X$. 
Then there exsists a positive integer $k_{0}$ depending only 
on $E$ such that for every $k\geqq k_{0}$
\[
H^{q}(X,{\cal O}_{X}(K_{X}+mL+kA)\otimes E) = 0
\]
holds for every $m\geqq 0$ and $q\geqq 1$. 
\end{lemma}
{\bf Proof}.
Let $\omega_{X}$ be the $L^{2}$-dualizing sheaf of 
$X$, i.e., the direct image sheaf of the canonical sheaf 
of a resolution of $X$. 
Since $X$ has only canonical singularities, we see that 
$\omega_{X}$ is isomorphic to ${\cal O}_{X}(K_{X})$. 
Since $L$ is nef and $A$ is ample, there exists a positive
integer $k_{0}$ such that for every $k\geqq k_{0}$, 
$(mL+kA)\otimes E$ admits a $C^{\infty}$-hermitian metric 
with (strictly) Nakano positive curvature. 
Then by the $L^{2}$-vanishing theorem, 
we see that 
\[
H^{q}(X,{\cal O}_{X}(K_{X}+mL+kA)\otimes E) = 0
\]
holds for every $m\geqq 0$ and $q\geqq 1$. 
This completes the proof. 
{\bf Q.E.D.}
 
Author's address\\
Hajime Tsuji\\
Department of Mathematics\\
Sophia University\\
7-1 Kioicho, Chiyoda-ku 102-8554\\
Japan \\
e-mail address: tsuji@mm.sophia.ac.jp


\begin{thebibliography}{99}
\bibitem{a-s} U. Anghern-Y.-T. Siu, Effective freeness and point separation
for adjoint bundles, Invent. Math. 122 (1995), 291-308.
\bibitem{b-d} T. Bandeman- G. Detholoff, Estimates on the number
of rational maps from a fixed variety to varieties of general type, Ann. Inst. Fourier 47, 
801-824 (1997). 
\bibitem{b3} E. Bombieri, Canonical models of surfaces of general type,
Publ. I.H.E.S. 42 (1972), 171-219.
\bibitem{dem} J.P. Demailly,  A numerical criterion for very ample line bundles, J. Diff. Geom. 37 (1993), 323-374.
\bibitem{d} J.P. Demailly,  Regularization of closed positive currents and 
intersection theory, J. of Alg. Geom. 1 (1992) 361-409.
\bibitem{d-p-s}  J.P. Demailly-T. Peternell-M. Schneider : 
Pseudo-effective line bundles on compact K\"{a}hler manifolds, 
math. AG/0006025 (2000).
\bibitem{f}A. R. Fletcher, Contributions to Riemann-Roch on projective $3$-folds with only canonical singularities and applications. Algebraic geometry, Bowdoin, 1985 (Brunswick, Maine, 1985), 221--231, Proc. Sympos. Pure Math., 46, Part 1, Amer. Math. Soc., Providence, RI, 1987. 
\bibitem{fu}T.  Fujita, On the structure of polarized varieties with $\Delta $-genera zero. J. Fac. Sci. Univ. Tokyo Sect. IA Math. 22 (1975), 103--115. 
\bibitem{ha} R. Harthshorne, Algebraic geometry, GTM 52, Springer 
(1977). 
\bibitem{h} L. H\"{o}rmander, An Introduction to Complex Analysis in Several
Variables 3-rd ed.,North-Holland(1990).
\bibitem{ka} Y. Kawamata, Subadjunction of log canonical divisors II, alg-geom math.AG/9712014, Amer. J. of Math. 120 (1998),893-899.
\bibitem{ka2} Y. Kawamata, Fujita's freeness conjecture for 3-folds and 
4-folds, Math. Ann. 308 (1997), 491-505.  
\bibitem{ka3} Y. Kawamata, Semipositivity, vanishing and applications, 
Lecture note in School ON VANISHING THEOREMS AND EFFECTIVE RESULTS IN 
ALGEBRAIC GEOMETRY (ICTP, Trieste, May 2000). 
\bibitem{k-o} S. Kobayashi-T. Ochiai, Mappings into compact complex manifolds
with negative first Chern class, Jour. Math. Soc. Japan 23 (1971),137-148.
\bibitem{k-m} J. Koll\'{a}r- S. Mori, Birational geometry of algebraic varieties, Cambridge Tracts in Math., Cambridge University Press (1998).
\bibitem{kr} S. Krantz, Function theory of several complex variables, 
John Wiley and Sons (1982).  
\bibitem{l}T. Luo,  Global $2$-forms on regular $3$-folds of general type. Duke Math. J. 71 (1993), no. 3, 859--869.
\bibitem{l2}T. Luo,Global holomorphic 2-forms and pluricanonical systems on threefolds. Math. Ann. 318 (2000), no. 4, 707--730. 
\bibitem{m} K. Maehara, A finiteness properties of varieties of general type,
Math Ann. 262 (1983),101-123.
\bibitem{mo} S. Mori, Flip conjecture and the existence of minimal
model for 3-folds, J. of A.M.S. 1 (1988), 117-253.
\bibitem{n}A.M. Nadel, Multiplier ideal sheaves and existence of K\"{a}hler-Einstein
metrics of positive scalar curvature, Ann. of Math. 132 (1990),549-596.
\bibitem{nak} N. Nakayama, Invariance of plurigenera of algebraic varieties,
 RIMS preprint 1191, March (1998). 
\bibitem{o-t}T. Ohsawa and K. Takegoshi, $L^{2}$-extension of holomorphic
functions, Math. Z. 195 (1987),197-204.
\bibitem{o} T. Ohsawa, On the extension of $L^{2}$ holomorphic functions V,
effects of generalization, Nagoya Math. J. (2001) 1-21. 
\bibitem{ti}G. Tian, On a set of polarized K\"{a}hler metrics on algebraic
manifolds, Jour. Diff. Geom. 32 (1990),99-130.
\bibitem{tu}H. Tsuji, Analytic Zariski decomposition, Proc. of Japan Acad.
61(1992) 161-163.
\bibitem{tu2} H. Tsuji, Existence and Applications of Analytic Zariski Decompositions, Trends in Math. Analysis and Geometry in Several Complex Variables,  (1999) 253-272.
\bibitem{tu3} H. Tsuji, Numeircally trivial fibrations, 
math.AG/0001023(2000). 
\bibitem{tu4} H. Tsuji, On the structure of pluricanonical systems of projective varieties of general type, preprint (1997).
\bibitem{t}H. Tsuji, Global generation of adjoint bundles, Nagoya Math. J.
142 (1996),5-16.
\bibitem{tu6} H. Tsuji, Deformation invariance of plurigenera, Nagoya Math. J. 166 (2002), 117-134. 
\bibitem{tu5} H. Tsuji, Subadjunction theorem for pluricanonical divisors, 
math.AG/0111311 (2001). 
\bibitem{3} H. Tsuji, Pluricanonical systems of projective 
 3-folds of general type, preprint (2002). 
\end{thebibliography}
\end{document}